\documentclass[11pt]{amsart}
\usepackage[letterpaper, total={6.5in,9in}, includefoot, centering]{geometry}

\usepackage{soul}
\usepackage{hyperref}
\usepackage{url}
\usepackage{amsmath,amssymb,amsfonts,stmaryrd,mathtools}
\usepackage[foot]{amsaddr}
\usepackage{bm}
\usepackage{xfrac}
\usepackage{array}
\usepackage{blkarray}
\usepackage{arydshln}
\usepackage{multirow}
\usepackage{empheq}
\usepackage{enumitem}
\usepackage{comment}

\usepackage{color}
\usepackage{adjustbox}
\definecolor{darkgreen}{rgb}{0.0, 0.5, 0}
\usepackage[numbers,sort]{natbib}
\usepackage{cleveref}
\usepackage[makeroom]{cancel}

\usepackage{xcolor}

\makeatother
\newtheorem{theorem}{Theorem}[section]

\newtheorem{definition}{Definition}[section]
\newtheorem{remark}{Remark}[section]
\newtheorem{example}{Example}[section]
\newtheorem{proposition}{Proposition}[section]

\graphicspath{ {./figures/} }

\usepackage{tikz}
\usetikzlibrary{fit}
\usetikzlibrary{shapes.geometric}

\usepackage{algorithm}
\usepackage{algpseudocode}

\usepackage{hyperref}
\hypersetup{colorlinks=true, urlcolor=blue, citecolor=blue}

\begin{document}
\def\treetype{2} 

\if\treetype1 

	\def\T#1#2{\scalebox{1}{\enspace\raisebox{-1mm}{\csname Tree#1#2\endcsname}\enspace}}
	\tikzstyle{nodeB}=[circle,fill=black,draw=black,scale=0.20,line width=0.8pt,radius=0.7pt]
	\tikzstyle{nodeW}=[circle,fill=black,draw=black,scale=0.24,line width=0.8pt,radius=0.7pt,fill=white]
	
	
	\tikzstyle{treeX}=[rotate=180,scale=0.15,line width=1.1pt]
	\expandafter\newcommand\csname Tree10\endcsname{\begin{tikzpicture}[treeX]\node[nodeB]{};\end{tikzpicture}}
	
	\expandafter\newcommand\csname Tree11\endcsname{\begin{tikzpicture}[treeX]\node[nodeW]{};\end{tikzpicture}}
	
	\expandafter\newcommand\csname Tree20\endcsname{\begin{tikzpicture}[treeX]\node[nodeB]{}child{node[nodeB]{}};\end{tikzpicture}}
	\expandafter\newcommand\csname Tree21\endcsname{\begin{tikzpicture}[treeX]\node[nodeB]{}child{node[nodeW]{}};\end{tikzpicture}}
	
	\expandafter\newcommand\csname Tree22\endcsname{\begin{tikzpicture}[treeX]\node[nodeW]{}child{node[nodeB]{}};\end{tikzpicture}}
	\expandafter\newcommand\csname Tree23\endcsname{\begin{tikzpicture}[treeX]\node[nodeW]{}child{node[nodeW]{}};\end{tikzpicture}}

	\expandafter\newcommand\csname Tree30\endcsname{\begin{tikzpicture}[treeX]\node[nodeB]{}child{node[nodeB]{}}child{node[nodeB]{}};\end{tikzpicture}}
	\expandafter\newcommand\csname Tree31\endcsname{\begin{tikzpicture}[treeX]\node[nodeB]{}child{node[nodeW]{}}child{node[nodeB]{}};\end{tikzpicture}}
	\expandafter\newcommand\csname Tree32\endcsname{\begin{tikzpicture}[treeX]\node[nodeB]{}child{node[nodeB]{}}child{node[nodeW]{}};\end{tikzpicture}}
	\expandafter\newcommand\csname Tree33\endcsname{\begin{tikzpicture}[treeX]\node[nodeB]{}child{node[nodeW]{}}child{node[nodeW]{}};\end{tikzpicture}}
	
	\expandafter\newcommand\csname Tree34\endcsname{\begin{tikzpicture}[treeX]\node[nodeW]{}child{node[nodeB]{}}child{node[nodeB]{}};\end{tikzpicture}}
	\expandafter\newcommand\csname Tree35\endcsname{\begin{tikzpicture}[treeX]\node[nodeW]{}child{node[nodeW]{}}child{node[nodeB]{}};\end{tikzpicture}}
	\expandafter\newcommand\csname Tree36\endcsname{\begin{tikzpicture}[treeX]\node[nodeW]{}child{node[nodeB]{}}child{node[nodeW]{}};\end{tikzpicture}}
	\expandafter\newcommand\csname Tree37\endcsname{\begin{tikzpicture}[treeX]\node[nodeW]{}child{node[nodeW]{}}child{node[nodeW]{}};\end{tikzpicture}}

	\expandafter\newcommand\csname Tree40\endcsname{\begin{tikzpicture}[treeX]\node[nodeB]{}child{node[nodeB]{}child{node[nodeB]{}}};\end{tikzpicture}}
	\expandafter\newcommand\csname Tree41\endcsname{\begin{tikzpicture}[treeX]\node[nodeB]{}child{node[nodeB]{}child{node[nodeW]{}}};\end{tikzpicture}}
	\expandafter\newcommand\csname Tree42\endcsname{\begin{tikzpicture}[treeX]\node[nodeB]{}child{node[nodeW]{}child{node[nodeB]{}}};\end{tikzpicture}}
	\expandafter\newcommand\csname Tree43\endcsname{\begin{tikzpicture}[treeX]\node[nodeB]{}child{node[nodeW]{}child{node[nodeW]{}}};\end{tikzpicture}}
	
	\expandafter\newcommand\csname Tree44\endcsname{\begin{tikzpicture}[treeX]\node[nodeW]{}child{node[nodeB]{}child{node[nodeB]{}}};\end{tikzpicture}}
	\expandafter\newcommand\csname Tree45\endcsname{\begin{tikzpicture}[treeX]\node[nodeW]{}child{node[nodeB]{}child{node[nodeW]{}}};\end{tikzpicture}}
	\expandafter\newcommand\csname Tree46\endcsname{\begin{tikzpicture}[treeX]\node[nodeW]{}child{node[nodeW]{}child{node[nodeB]{}}};\end{tikzpicture}}
	\expandafter\newcommand\csname Tree47\endcsname{\begin{tikzpicture}[treeX]\node[nodeW]{}child{node[nodeW]{}child{node[nodeW]{}}};\end{tikzpicture}}

\fi

\if\treetype2 

	
	
	\def\T#1#2{\scalebox{1.5}{\enspace\raisebox{-1mm}{\csname Tree#1#2\endcsname}\enspace}}
	\tikzstyle{dns}=[circle,fill=black,draw=black,scale=0.20,line width=0.5pt,solid] 
	\tikzstyle{edgeB}=[black,line width=1pt,solid]
	\tikzstyle{edgeW}=[black,line width=0.9pt,densely dotted]
	
	
	\tikzstyle{treeX}=[rotate=180,scale=0.15,line width=1.1pt]
	\expandafter\newcommand\csname Tree10\endcsname{\begin{tikzpicture}[treeX]\node[dns]{};\end{tikzpicture}}
	
	\expandafter\newcommand\csname Tree11\endcsname{\begin{tikzpicture}[treeX]\node[dns]{};\end{tikzpicture}}
	
	\expandafter\newcommand\csname Tree20\endcsname{\begin{tikzpicture}[treeX]\node[dns]{}child[edgeB]{node[dns]{}};\end{tikzpicture}}
	\expandafter\newcommand\csname Tree21\endcsname{\begin{tikzpicture}[treeX]\node[dns]{}child[edgeW]{node[dns]{}};\end{tikzpicture}}
	
	\expandafter\newcommand\csname Tree22\endcsname{\begin{tikzpicture}[treeX]\node[dns]{}child{node[dns]{}};\end{tikzpicture}}
	\expandafter\newcommand\csname Tree23\endcsname{\begin{tikzpicture}[treeX]\node[dns]{}child{node[dns]{}};\end{tikzpicture}}

	\expandafter\newcommand\csname Tree30\endcsname{\begin{tikzpicture}[treeX]\node[dns]{}child{node[dns]{}}child{node[dns]{}};\end{tikzpicture}}
	\expandafter\newcommand\csname Tree31\endcsname{\begin{tikzpicture}[treeX]\node[dns]{}child[edgeW]{node[dns]{}}child{node[dns]{}};\end{tikzpicture}}
	\expandafter\newcommand\csname Tree32\endcsname{\begin{tikzpicture}[treeX]\node[dns]{}child{node[dns]{}}child[edgeW]{node[dns]{}};\end{tikzpicture}}
	\expandafter\newcommand\csname Tree33\endcsname{\begin{tikzpicture}[treeX]\node[dns]{}child[edgeW]{node[dns]{}}child[edgeW]{node[dns]{}};\end{tikzpicture}}
	
	\expandafter\newcommand\csname Tree34\endcsname{\begin{tikzpicture}[treeX]\node[dns]{}child{node[dns]{}}child{node[dns]{}};\end{tikzpicture}}
	\expandafter\newcommand\csname Tree35\endcsname{\begin{tikzpicture}[treeX]\node[dns]{}child{node[dns]{}}child{node[dns]{}};\end{tikzpicture}}
	\expandafter\newcommand\csname Tree36\endcsname{\begin{tikzpicture}[treeX]\node[dns]{}child{node[dns]{}}child{node[dns]{}};\end{tikzpicture}}
	\expandafter\newcommand\csname Tree37\endcsname{\begin{tikzpicture}[treeX]\node[dns]{}child{node[dns]{}}child{node[dns]{}};\end{tikzpicture}}

	\expandafter\newcommand\csname Tree40\endcsname{\begin{tikzpicture}[treeX]\node[dns]{}child{node[dns]{}child{node[dns]{}}};\end{tikzpicture}}
	\expandafter\newcommand\csname Tree41\endcsname{\begin{tikzpicture}[treeX]\node[dns]{}child{node[dns]{}child[edgeW]{node[dns]{}}};\end{tikzpicture}}
	\expandafter\newcommand\csname Tree42\endcsname{\begin{tikzpicture}[treeX]\node[dns]{}child[edgeW]{node[dns]{}child[edgeB]{node[dns]{}}};\end{tikzpicture}}
	\expandafter\newcommand\csname Tree43\endcsname{\begin{tikzpicture}[treeX]\node[dns]{}child[edgeW]{node[dns]{}child[edgeW]{node[dns]{}}};\end{tikzpicture}}
	
	\expandafter\newcommand\csname Tree44\endcsname{\begin{tikzpicture}[treeX]\node[dns]{}child{node[dns]{}child{node[dns]{}}};\end{tikzpicture}}
	\expandafter\newcommand\csname Tree45\endcsname{\begin{tikzpicture}[treeX]\node[dns]{}child{node[dns]{}child{node[dns]{}}};\end{tikzpicture}}
	\expandafter\newcommand\csname Tree46\endcsname{\begin{tikzpicture}[treeX]\node[dns]{}child{node[dns]{}child{node[dns]{}}};\end{tikzpicture}}
	\expandafter\newcommand\csname Tree47\endcsname{\begin{tikzpicture}[treeX]\node[dns]{}child{node[dns]{}child{node[dns]{}}};\end{tikzpicture}}

\fi
\title{Order Conditions for Nonlinearly Partitioned Runge--Kutta Methods}
\author{Brian K. Tran$^\dagger$ and Ben S. Southworth$^\dagger$ and Tommaso Buvoli$^\ddagger$}
\address{$^\dagger$Los Alamos National Laboratory, Theoretical Division, Los Alamos, NM 87545, USA. \newline \indent $^\ddagger$Tulane University, Department of Mathematics, \newline \indent\ \ New Orleans, LA 70118, USA.}
\email{btran@lanl.gov, southworth@lanl.gov, tbuvoli@tulane.edu}
\allowdisplaybreaks

\begin{abstract}
Recently a new class of \emph{nonlinearly partitioned Runge--Kutta} (NPRK) methods was proposed for nonlinearly partitioned systems of autonomous ordinary differential equations, $y' = F(y,y)$. The target class of problems are ones in which different scales, stiffnesses, or physics are coupled in a nonlinear way, wherein the desired partition cannot be written in a classical additive or component-wise fashion. Here we use rooted-tree analysis to derive full order conditions for NPRK$_M$ methods, where $M$ denotes the number of nonlinear partitions. Due to the nonlinear coupling and thereby mixed product differentials, it turns out the standard node-colored rooted-tree analysis used in analyzing ODE integrators does not naturally apply. Instead we develop a new edge-colored rooted-tree framework to address the nonlinear coupling. The resulting order conditions are enumerated, provided directly for up to 4th order with $M=2$ and 3rd-order with $M=3$, and related to existing order conditions of additive and partitioned RK methods. We conclude with an example which shows how the nonlinear order conditions can be used to obtain an embedded estimate of the state-dependent nonlinear coupling strength in a dynamical system.
\end{abstract}

\maketitle

{
  \hypersetup{linkcolor=black}
  \tableofcontents
}

\section{Introduction}

Nonlinearly partitioned Runge-Kutta (NPRK) methods \cite{nprk1} are a newly-proposed family of time integrators for solving the \textit{partitioned} initial value problem
\begin{align}
	y' = F(y,y), \quad y(t_0) = y_0.
	\label{eq:ivp-nonlinearly-partitioned}
\end{align}
NPRK methods can treat each argument of $F(y,y)$ with a different level of implicitness, and may be interpreted as nonlinear generalizations of additive Runge-Kutta (ARK) methods \cite{Ascher.1997,Kennedy.2003tv4,Kennedy.2019}. An $s$-stage NPRK method for the partitioned initial value problem \eqref{eq:ivp-nonlinearly-partitioned} is
\begin{align}
	\begin{aligned}
		Y_i &= y_n + h\sum_{j=1}^{s} \sum_{k=1}^{s} a_{ijk} F(Y_j,Y_k), \quad i=1,\ldots, s, \\
		y_{n+1} &= y_n + h\sum_{i=1}^s \sum_{j=1}^s b_{ij} F(Y_i,Y_j),	
	\end{aligned}
	\label{eq:nprk-general}
\end{align}
where $Y_i$ are stage values, the rank three tensor $a_{ijk}$ takes the place of the classical Runge-Kutta matrix $a_{ij}$, and the matrix $b_{ij}$ replaces the classical weight vector $b_i$. To motivate the definition of the partitioned initial value problem \eqref{eq:ivp-nonlinearly-partitioned}, we first provide a simple motivating example (for another example, see \Cref{sec:example} where we consider the Lotka--Volterra equation).
\begin{example}\label{example:burger:example}[Viscous Burgers' Equation]
\textit{
    To provide intuition for how partitioned initial value problems can lead to more efficient methods, we consider the viscous Burgers' equation with viscosity coefficient $\epsilon \geq 0$ given by
    $$ u_t = \epsilon u_{xx} + u u_x, $$
    which, after spatial semi-discretization, takes the form
    \begin{equation}\label{eq:example-burger}
        \mathbf{u}' = \epsilon D \mathbf{u} + \textup{diag}(\mathbf{u}) A \mathbf{u} =: f(\mathbf{u}),
    \end{equation}
    where $D$ and $A$ denote the discrete diffusion and advection matrices, respectively. It is instructive to consider three standard first-order methods for integrating \eqref{eq:example-burger}.
    \begin{align*}
        \textup{Implicit: }& \mathbf{u}_{n+1} = \mathbf{u}_n + h f(\mathbf{u}_{n+1}), \\
        \textup{Rosenbrock: }& (I - hJ_n)(\mathbf{u}_{n+1} - \mathbf{u}_n) = hf(\mathbf{u}_n), \\
        \textup{IMEX: }& \mathbf{u}_{n+1} = \mathbf{u}_n + h \epsilon D \mathbf{u}_{n+1} + h \textup{diag}(\mathbf{u}_n) A \mathbf{u}_n,
    \end{align*}
    where here $J_n$ denotes the Jacobian of $f$ at $\mathbf{u}_n$. The implicit method is of course stable independent of $\epsilon$, but requires a nonlinear solve. The stability and accuracy of the Rosenbrock method are comparable to that of the implicit method, but only requires a linear solve involving the Jacobian \cite{nprk1}. The IMEX method requires only a linear solve; however, it reduces to explicit Euler when the viscosity coefficient $\epsilon \rightarrow 0$ and therefore, can suffer from stability issues for small $\epsilon$. In \cite{nprk1}, we show how formulating \eqref{eq:example-burger} as a partitioned initial value problem and subsequently applying an NPRK method can retain the advantages of each of the above methods while eliminating their respective disadvantages. Namely, note that \eqref{eq:example-burger} can be expressed as a partitioned initial value problem $\mathbf{u}' = F(\mathbf{u},\mathbf{u})$ where
    $$ F(\mathbf{u},\mathbf{v}) := \epsilon D \mathbf{u} + \textup{diag}(\mathbf{v}) A \mathbf{u}.$$
    A corresponding first-order NPRK method is given by
    \begin{align*}
        \mathbf{u}_{n+1} &= \mathbf{u}_n + h F(\mathbf{u}_{n+1}, \mathbf{u}_n) \\
        &= \mathbf{u}_n + h \epsilon D \mathbf{u}_{n+1} + \textup{diag}(\mathbf{u}_{n}) A \mathbf{u}_{n+1}.
    \end{align*}
    Note that this method is linearly implicit for all $\epsilon$, provides the same accuracy and stability as the implicit and Rosenbrock method, while only requiring a linear solve at each step and avoiding computing the Jacobian (for numerical results comparing these methods, see \cite{nprk1}). 
}
\end{example}

The NPRK framework facilitates efficient integration of equations with stiff terms that cannot be isolated additively or component-wise. Specifically, any unpartitioned system $y'=f(y)$ can be converted into \eqref{eq:ivp-nonlinearly-partitioned} by selecting a function $F(y,z)$ that satisfies $f(y) = F(y,y)$; if integrated with an IMEX-NPRK method (i.e. $a_{ijk} = 0$ for $j>i$ and $k \ge i$), then implicit solves are only required over the first argument of $F$, e.g., as discussed in Example \ref{example:burger:example}. This approach has been previously applied to solve a number of problems using the related, though less general family of semi-implicit integrators \cite{Boscarino.2022,Boscheri.2022,southworth2023implicit}. In \cite{nprk1}, we had success applying NPRK methods to nonlinear partitions of challenging thermal radiative transfer and radiation hydrodynamics problems as studied in \cite{southworth2023implicit,imex-trt}, and have ongoing work applying the framework to other multiphysics systems.

This paper derives order conditions for the NPRK method \eqref{eq:nprk-general} and the more general NPRK$_M$ method, described later in \eqref{subsec:m-nprk-intro}, that allows for $M$ nonlinear partitions. We derive the order conditions using edge-colored rooted trees that represent the elementary differentials in the Taylor series expansion of \eqref{eq:ivp-nonlinearly-partitioned}. This work generalizes our previous results from \cite{nprk1} which only investigated order conditions for a simple subclass of NPRK methods called sequentially-coupled methods.

\subsection{Main contribution}

This paper is the first to analyze the full order conditions for the NPRK method \eqref{eq:nprk-general} and the more general class of NPRK$_M$ method \eqref{eq:nprk-general-m-partitions}. In our the previous work \cite{nprk1}, we only considered a simplified family of NPRK methods, whose order conditions are equivalent to those of ARK methods. In this work, we show that ARK order conditions are insufficient for the general class of NPRK methods, due to nonlinear coupling terms between the arguments of a nonlinear partition.

We first review NPRK methods in \Cref{section:nprk-review}. To derive order conditions for NPRK methods, we introduce a new edge-colored rooted tree framework, and use it to first derive order conditions for two partitions in \Cref{subsection: tree order conditions}, Theorem \ref{thm:2-nprk-order-cond}, and generalize this to $M$-partitions in \Cref{sec:order:n}, Theorem \ref{thm:n-nprk-order-cond}. Full order conditions for up to five partitions and 8th order are enumerated in \Cref{sec:order:enum}. Interestingly, an NPRK$_M$ method has $1/M$ times the number of order conditions of a $\text{ARK}_M$  method, which is due to the inherent tensorial structure of an NPRK tableau. Furthermore, in \Cref{sec:order:enum}, we present an algorithm for computing NPRK order conditions, and explicit order conditions up to 4th order for $M=2$ and 3rd order for $M=3$ are provided in \Cref{app:order-list}.

NPRK methods and the underlying order conditions are then related to ARK and PRK methods in \Cref{section:relation-to-ark-order-cond}. By relating NPRK order conditions to ARK order conditions, we show that NPRK order conditions contain new order conditions corresponding to nonlinear coupling which vanish when considering additive partitions. As an application of these nonlinear order conditions, in \Cref{sec:example}, we show through an example how to obtain an embedded estimate for the nonlinear coupling strength in a system. The code for all of the numerical results obtained in the paper, including numerical and symbolic implementations of the NPRK order conditions, and for the numerical example, is available at \cite{githubNPRK2025}.

\section{Nonlinearly partitioned Runge--Kutta methods}\label{section:nprk-review}

We introduce several definitions and properties of NPRK methods \eqref{eq:nprk-general} that are relevant for determining and analyzing order conditions.

\subsection{Nonlinear partitions}

Consider an ordinary differential equation (ODE) $\dot{y} = f(y)$ on a vector space $X$, specified by a vector field $f: X \rightarrow X$. We say that a mapping $F: X \times X \rightarrow X$ is a \emph{nonlinear partition} of $f: X \rightarrow X$ if
    \begin{align}
        F(y,y) = f(y) \text{ for all } y \in X.
        \label{eq:partition}
    \end{align}
Nonlinear partitions for a given $f(y)$ are not unique, and different choices will affect the stability, accuracy, and computational efficiency of an NPRK method. For example, even a simple scalar function like $f(y) = y^2$ can be nonlinearly partitioned in an infinite number of ways; $F(u,v) = u^p v^{2-p}$, $p \in \mathbb{R}$ is one such example.

\subsection{Underlying RK and ARK integration methods}\label{sec:underlying-rk-ark}

For certain nonlinear partitions, NPRK methods \eqref{eq:nprk-general} reduce to simpler RK families.  
\begin{definition}
	A classical RK integrator $(a_{ij}, b_i, c_i)$ is an \emph{underlying RK method} of an NPRK integrator, if \eqref{eq:nprk-general} reduces to the Runge--Kutta method corresponding to $(a_{ij}, b_i, c_i)$ when the function $F(x,y)$ depends only on $x$ or only on $y$.
\end{definition}

If ${F(u,v) = f(u)}$, then \eqref{eq:nprk-general} reduces to the classical RK method
\begin{align}
        a^{\{1\}}_{ij} = \sum_{k=1}^s a_{ijk}, \quad  
        b^{\{1\}}_{i} = \sum_{j=1}^s b_{ij}, \quad
        c^{\{1\}}_i = \sum_{j=1}^{s} a^{\{1\}}_{ij}.
    \label{eq:nprk-underlying-rk1}
\end{align}
Similarly, if ${F(u,v) = f(v)}$, then \eqref{eq:nprk-general} reduces to
\begin{align}
        a^{\{2\}}_{ik} = \sum_{j=1}^s a_{ijk}, \quad
        b^{\{2\}}_{j} = \sum_{i=1}^s b_{ij}, \quad
        c^{\{2\}}_i = \sum_{k=1}^{s} a^{\{2\}}_{ik} = c^{\{1\}}_i.
    \label{eq:nprk-underlying-rk2}
\end{align}
If the NPRK method \eqref{eq:nprk-general} is applied to a nonlinear partition with additive structure
    \begin{align}
    F(u,v)=F^{\{1\}}(u) +F^{\{2\}}(v);
    \label{eq:nprk-underlying-ark-equation}
\end{align} 
then it reduces to an ARK method \cite{Kennedy.2003tv4}, which we call the {\em underlying ARK method}.
\begin{definition}
	The underlying ARK method of \eqref{eq:nprk-general}	is
	\begin{align}
   \begin{aligned}
   Y_{i} &= y_n + h \sum_{j=1}^s \left[ a^{\{1\}}_{ij} F^{\{1\}}(Y_j) + a^{\{2\}}_{ij} F^{\{2\}}(Y_j) \right], \\
   y_{n+1} &= y_n + h \sum_{i=1}^s \left[ b^{\{1\}}_{i} F^{\{1\}}(Y_i) + b^{\{2\}}_{i} F^{\{2\}}(Y_i) \right],
  \end{aligned}
  \label{eq:nprk-underlying-ark}
\end{align}
whose coefficients are defined in \eqref{eq:nprk-underlying-rk1} and \eqref{eq:nprk-underlying-rk2}. 
\end{definition}

Since NPRK methods maintain their order-of-accuracy for all nonlinear partitions, it follows that an order $q$ accurate NPRK method, must have underlying methods (RK and ARK) that are also order $q$ accurate. In \Cref{sec:order}, we will distinguish the resulting RK and ARK order conditions from those that are due to nonlinear coupling between the arguments of a nonlinear partition.

\subsection{Generality of ARK and NPRK}

NPRK methods generalize ARK methods in the sense that they allow for nonlinearly partitioned right-hand-sides \eqref{eq:ivp-nonlinearly-partitioned}, but reduce to their underlying ARK method \eqref{eq:nprk-underlying-ark} on additively partitioned right-hand-sides. We briefly discuss the generality of these two method families.

	\begin{itemize}
		\item {\em NPRK to ARK.}	 The NPRK family only contains ARK methods with shared abscissa $c^{\{1\}}_i = c^{\{2\}}_i$. An ARK method with differing abscissa cannot be an underlying ARK method because the abscissa of all underlying ARK methods are all formed by summing over the same indices of the tensor $a_{ijk}$:
			\begin{align*}
				\textstyle c_i^{\{1\}}=\sum_j \big( \sum_k a_{ijk} \big) = \sum_k \big( \sum_j a_{ijk} \big) = c_i^{\{2\}}. 
			\end{align*}
		\sloppy \item {\em ARK to NPRK.} Given any ARK method $\mathcal{M}$ with tableaux $(a^{\{1\}},b^{\{1\}},c^{\{1\}})$, $(a^{\{2\}},b^{\{2\}},c^{\{2\}})$, shared abscissa  $c^{\{1\}}_i = c^{\{2\}}_i$ and $s>1$ stages, it is possible to construct a family of NPRK methods with $\mathcal{M}$ as their underlying ARK method. This restriction imposes $2s^2$ conditions on the tensor $a_{ijk}$ and $2s$ conditions on the matrix $b_j$, leading to a total of $(s^3 - 2s^2)$ free parameters in $a_{ijk}$ and $s^2 - 2s$ free parameters in $b_{ij}$.
	\end{itemize}

\subsection{$M$-nonlinearly partitioned Runge-Kutta}
\label{subsec:m-nprk-intro}
The NPRK framework trivially generalizes to the $M$-component nonlinearly partitioned equation $y'=F(y,\ldots, y)$ where $F$ has $M$ arguments. The ansatz for an $s$-stage, NPRK$_M$ method is
\begin{align}
	\begin{aligned}
		Y_i &= y_n + h\sum_{j_1,\ldots, j_M=1}^{s} a_{i, j_1,\ldots, j_M} F(Y_{j_1}, \ldots, Y_{j_m}), \quad i=1,\ldots, s, \\
		y_{n+1} &= y_n + h \sum_{j_1,\ldots, j_M=1}^{s} b_{j_1, \ldots, j_M} F(Y_{j_1}, \ldots, Y_{j_M}),	
	\end{aligned}
	\label{eq:nprk-general-m-partitions}
\end{align}
where the tensors $a$ and $b$ are now of rank $M+1$ and $M$, respectively. The method \eqref{eq:nprk-general-m-partitions} nonlinearly combines $M$ underlying RK methods and reduces to an $\text{ARK}_M$ method \cite{Kennedy.2003tv4} when applied to an additively partitioned system $F(Y_{1}, \ldots, Y_{M}) = \sum_{k=1}^M F^{\{k\}}(Y_{k})$.

\section{NPRK order conditions}\label{sec:order}

In this section, we present order conditions for the NPRK method \eqref{eq:nprk-general} and the more general NPRK$_M$ method \eqref{eq:nprk-general-m-partitions}. We begin by computing third-order NPRK conditions via Taylor expansion, and then introduce edge-color rooted trees for determining order conditions for NPRK$_M$ methods. 

\subsection{Order conditions from Taylor series}\label{sec:order:taylor}

The most direct approach for deriving order conditions is to compare successive terms in the Taylor expansion of the exact solution to those of an NPRK method. In doing so, we will also see that the derivatives, and order conditions are related to rooted trees with colored edges.

Let $y(t)$ denote the exact solution to the initial-value problem \eqref{eq:ivp-nonlinearly-partitioned}. To Taylor expand $y(t_0 + h)$ about $t_0$, we first compute the derivatives of $y$ up to third order,
\begin{align*}
    y' &= F(y,y), \\
    y'' &= D_1F(y,y)y' + D_2F(y,y)y', \\
    y''' &= D_{11} F(y,y)[y',y'] + 2D_{12} F(y,y)[y',y'] + D_{22} F(y,y)[y',y'] + D_1F(y,y) D_1F(y,y) y'\\ 
    & \qquad  + D_1F(y,y)D_2F(y,y)y' + D_2F(y,y)D_1F(y,y)y' + D_2F(y,y)D_2F(y,y)y'.
\end{align*}
Here, $D_iF(y_1,y_2): X \rightarrow X$ denotes the derivative of $F$ at $(y_1,y_2)$ with respect to the $i^{th}$ argument, which is a linear mapping whose action on $z \in X$ we simply write as $D_iF(y_1,y_2)z.$ Similarly, $D_{ij}F(y_1,y_2): X \times X \rightarrow X$ denotes the second derivative of $F$ at $(y_1,y_2)$ with respect to the $j^{th}$ then $i^{th}$ argument, which is a bilinear mapping, whose action on $u,v \in X$ we write as $D_{ij}F(y_1,y_2)[u,v].$ Note, in the above, any instance of $y'$ appearing on the right hand side can be replaced with $F(y,y)$.

For brevity, throughout this section, unless the argument of $F$ or its derivatives are written explicitly, we will assume that they are evaluated at $(y_0,y_0)$. Thus, the Taylor expansion of $y(t_0+h)$ about $t_0$ is given by
\begin{align} \label{eq:third-order-taylor-series}
    y(t_0+h) &= y(t_0) + hy'(t_0) + \frac{h^2}{2}y''(t_0) + \frac{h^2}{3!}y'''(t_0) + \mathcal{O}(h^4) \\
    &= y_0 + h F + \frac{h^2}{2} \Big(D_1F F + D_2F F \Big) \nonumber \\
    & \qquad + \frac{h^3}{3!}\Big[ D_{11} F[F,F] + 2D_{12} F[F,F] + D_{22} F[F,F] + D_1F D_1F F \nonumber\\ 
    & \qquad \qquad \qquad + D_1F D_2F F + D_2FD_1FF + D_2FD_2FF \Big] + \mathcal{O}(h^4). \nonumber
\end{align}

Each summand in the derivatives (e.g., $D_1 FF$ or $D_{12}F[F,F]$) is called an {\em elementary differential}. We can find a pattern in the elementary differentials using edge-colored rooted trees, as shown in \Cref{tab:edge-colored-rooted-trees-intro}; we formalize this relationship in \Cref{subsection: tree order conditions}. For the moment, we simply remark that: (1) each factor in an elementary differential corresponds to a node, (2) the derivative order ($1$ for $D_i$ and $2$ for $D_{ij}$) determines the number of outward edges, and (3) the edge color depends on whether differentiation is conducted with respect to the first or second argument. 

\begin{table}[H]
	\begin{scriptsize}
		\begin{tabular}{@{}c@{}}
			\begin{tabular}{|c|c|} \hline	
				Order 1 & Order 2 \\ \hline
				& \\
				\begin{tabular}{c}
					$\T10$ \\[1em]
					$F$
				\end{tabular} &
				\begin{tabular}{cc}
					$\T21$ 	& $\T20$ \\[1em]
					$F_1[F]$	& $F_2[F]$ \\[-0.5em] \\
				\end{tabular} \\ \hline
			\end{tabular} \\ \\ \hline
			\begin{tabular}{|c|}
				Order 3	\\\hline \\
				\setlength{\tabcolsep}{5.4pt}
				\begin{tabular}{@{}cccccccc@{}}
					$\T33$			& $\T32$ 		& $\T31$ 		& $\T30$  		& $\T43$			& $\T42$			& $\T41$			& $\T40$ \\[1em]
					$F_{11}[F,F]$	& $F_{12}[F,F]$	& $F_{21}[F,F]$	& $F_{22}[F,F]$	& $F_{1}[F_{1}[F]]$	& $F_{1}[F_{2}[F]]$	& $F_{2}[F_{1}[F]]$ & $F_{2}[F_{2}[F]]$ \\[-0.5em] \\
				\end{tabular}
			\end{tabular}  \\ \hline 
		\end{tabular}
	\end{scriptsize}
	\vspace{1em}

	\caption{Derivatives of $y$ up to order three and associated edge-colored rooted trees. The notation $F_{i}$ and $F_{ij}$ respectively abbreviates $D_i F$ and $D_{ij} F$. Thin dashed edges represent differentiation with respect to the first argument, while thick solid edges represent differentiation with respect to the second argument. Note that the third-order trees for $F_{12}[F,F]$ and $F_{21}[F,F]$ lead to a single order condition.}	
	\label{tab:edge-colored-rooted-trees-intro}

\end{table}

Next we expand the numerical solution obtained by applying one step of \eqref{eq:nprk-general}; a direct calculation yields
\begin{align}\label{eq:taylor-exp}
    y_1 & = y_0 + h \sum_{ij}b_{ij} F + h^2 \sum_{ijkl}b_{ij}a_{ikl}  D_1F F + h^2 \sum_{ijkl}b_{ij} a_{jkl} D_2F F \nonumber \\ 
    & \qquad + h^3 \sum_{ijkluv} b_{ij} a_{ikl} a_{kuv} D_1F D_1F F + h^3 \sum_{ijkluv} b_{ij} a_{ikl} a_{luv} D_1F D_2F F \nonumber \\ 
    & \qquad + h^3 \sum_{ijkluv} b_{ij} a_{jkl} a_{kuv} D_2F D_1F F + h^3 \sum_{ijkluv} b_{ij} a_{jkl} a_{luv} D_2F D_2F F \nonumber \\
    & \qquad + \frac{h^3}{2} \sum_{ijkluv}b_{ij}a_{ikl}a_{iuv} D_{11}F [F,F] + \frac{h^3}{2} \sum_{ijkluv} b_{ij}a_{ikl}a_{juv} 2 D_{12}F [F,F] \nonumber\\
    & \qquad + \frac{h^3}{2} \sum_{ijkluv} b_{ij}a_{jkl}a_{juv} D_{22}F[F,F] + \mathcal{O}(h^4). \nonumber
\end{align}
Comparing the above expansions of $y(t_0+h)$ and $y_1$ yields the order conditions up to order three in \Cref{tab:nprk-order-3-v2}.

\begin{table}[H]

	\begin{center}
		\renewcommand*{\arraystretch}{1.5}	
		\begin{tabular}{c|c|ll|l}
			Tree 	& Elementary 	& Elementary 			& $\Phi(\tau)$ in terms of 								& $\gamma(\tau)$ \\[-0.5em]
			$\tau$ 	& Differential 	& Weight $\Phi(\tau)$ 	& \eqref{eq:nprk-underlying-rk1}, \eqref{eq:nprk-underlying-rk2}  &  \\[0.5em] \hline
			$\T10$ & $F$ 				& $\sum_{ij}b_{ij}$ 						& $= \sum_i b_i^{\{1\}} = \sum_{i} b_i^{\{2\}}$ & 1 \\ \hdashline
			$\T21$ & $F_1[F]$  			& $\sum_{ijkl}b_{ij}a_{ikl}$ 			& $= \sum_{i} b_i^{\{1\}}c_i$  & 2 \\
			$\T20$ & $F_2[F]$  			& $\sum_{ijkl}b_{ij}a_{jkl}$ 			& $= \sum_{i} b_i^{\{2\}}c_i$ & 2 \\ \hdashline
			$\T33$ & $F_{11}[F,F]$ 		& $\sum_{ijkluv} b_{ij} a_{ikl} a_{iuv}$ & $= \sum_{i} b_i^{\{1\}} c_i c_i$ & 3 \\
			$\T32$ & $F_{12}[F,F]$ 		& $\sum_{ijkluv} b_{ij} a_{ikl} a_{juv}$ & n.a. & 3\\
			$\T30$ & $F_{22}[F,F]$ 		& $\sum_{ijkluv} b_{ij} a_{jkl} a_{juv}$ 	& $= \sum_{i} b_i^{\{2\}} c_i c_i$ & 3 \\[1em]
			$\T43$ & $F_{1}[F_{1}[F]]$ 	& $\sum_{ijkluv} b_{ij} a_{ikl} a_{kuv}$ 	& $= \sum_{ij} b_i^{\{1\}} a^{\{1\}}_{ij} c_j$ & 6\\[1em]
			$\T42$ & $F_{1}[F_{2}[F]]$ 	& $\sum_{ijkluv} b_{ij} a_{ikl} a_{luv}$ 	& $= \sum_{ij} b_i^{\{1\}} a^{\{2\}}_{ij} c_j$ & 6\\[1em]
			$\T41$ & $F_{2}[F_{1}[F]]$ 	& $\sum_{ijkluv} b_{ij} a_{jkl} a_{kuv}$ 	& $= \sum_{ij} b_i^{\{2\}} a^{\{2\}}_{ij} c_j$ & 6\\[1em]
			$\T40$ & $F_{2}[F_{2}[F]]$ 	& $\sum_{ijkluv} b_{ij} a_{jkl} a_{luv}$ 	& $= \sum_{ij} b_i^{\{2\}} a^{\{2\}}_{ij} c_j$ & 6\\
		\end{tabular}
		\vspace{1em}
	\end{center}

	\caption{Order conditions are $\Phi(\tau) = 1/\gamma(\tau)$. The horizontal dashed lines separate conditions of orders one, two, and three. 
    The order conditions consist of the the well-known ARK order conditions, along with two additional conditions corresponding to nonlinear coupling between the arguments of $F$.}
		\label{tab:nprk-order-3-v2}
	
\end{table}

\subsection{Order conditions via edge-colored rooted trees}\label{subsection: tree order conditions}

The starting point for a systematic derivation of NPRK order conditions is to observe that for $y' = F(y,y)$, the total time derivative is given by
$$ \frac{d}{dt}(\ ) = D_1(\ )\cdot F + D_2(\ )\cdot F, $$
where $\cdot F$ denotes insertion of $F$ into the multilinear operator obtained after differentiation.

Based on this observation, we can represent the derivatives of $y$ via edge-colored rooted trees. Every edge in a tree has one of two colors, which we represent graphically as dashed or solid, corresponding to differentiation with respect to the first or second argument, respectively, i.e., 
\begin{align*}
D_1(\ ) = \begin{tikzpicture}[baseline={0}]
    \draw[black, dashed, ultra thick] (0,0.4) -- (0,-0.4);
\end{tikzpicture}\ ,\ D_2(\ ) = \begin{tikzpicture}[baseline={0}]
    \draw[black, ultra thick] (0,0.4) -- (0,-0.4);
\end{tikzpicture}\ .
\end{align*}
Letting $F = \begin{tikzpicture} \filldraw[black] (0,0) circle (2pt); \end{tikzpicture}$, we can differentiate by grafting onto a rooted tree a dashed edge connected to a node and a solid edge connected to a node,
\begin{align*}
D_1(\ )\cdot F = \begin{tikzpicture}[baseline={0}]
    \filldraw[black] (0,0.4) circle (2pt);
    \draw[black, dashed, ultra thick] (0,0.4) -- (0,-0.4);
\end{tikzpicture}\ ,\ D_2(\ )\cdot F = \begin{tikzpicture}[baseline={0}]
    \filldraw[black] (0,0.4) circle (2pt);
    \draw[black, ultra thick] (0,0.4) -- (0,-0.4);
\end{tikzpicture}\ .
\end{align*}
For example, the derivatives of $y$ up to third-order are given by
\begin{align*}
y' &= F = \begin{tikzpicture}
    \filldraw[black] (0,0) circle (2pt);
\end{tikzpicture}\ , \\
y'' &= D_1F F + D_2F F = \begin{tikzpicture}[baseline={0}]
    \filldraw[black] (0,0.4) circle (2pt);
    \filldraw[black] (0,-0.4) circle (2pt);
    \draw[black, dashed, ultra thick] (0,0.4) -- (0,-0.4);
\end{tikzpicture} + \begin{tikzpicture}[baseline={0}]
    \filldraw[black] (0,0.4) circle (2pt);
    \filldraw[black] (0,-0.4) circle (2pt);
    \draw[black, ultra thick] (0,0.4) -- (0,-0.4);
\end{tikzpicture}\ , \\
y''' &= \underbrace{\begin{tikzpicture}[baseline={0}]
    \filldraw[black] (0,0.8) circle (2pt);
    \filldraw[black] (0,0) circle (2pt);
    \filldraw[black] (0,-0.8) circle (2pt);
    \draw[black, dashed, ultra thick] (0,0) -- (0,-0.8);
    \draw[black, dashed, ultra thick] (0,0) -- (0,0.8);
\end{tikzpicture}}_{=D_1F D_1F F} + \underbrace{\begin{tikzpicture}[baseline={0}]
    \filldraw[black] (0,0.8) circle (2pt);
    \filldraw[black] (0,0) circle (2pt);
    \filldraw[black] (0,-0.8) circle (2pt);
    \draw[black, ultra thick] (0,0) -- (0,-0.8);
    \draw[black, dashed, ultra thick] (0,0) -- (0,0.8);
\end{tikzpicture}}_{=D_2FD_1FF} + \underbrace{\begin{tikzpicture}[baseline={0}]
    \filldraw[black] (0,0.8) circle (2pt);
    \filldraw[black] (0,0) circle (2pt);
    \filldraw[black] (0,-0.8) circle (2pt);
    \draw[black, dashed, ultra thick] (0,0) -- (0,-0.8);
    \draw[black, ultra thick] (0,0) -- (0,0.8);
\end{tikzpicture}}_{=D_1FD_2FF} + \underbrace{\begin{tikzpicture}[baseline={0}]
    \filldraw[black] (0,0.8) circle (2pt);
    \filldraw[black] (0,0) circle (2pt);
    \filldraw[black] (0,-0.8) circle (2pt);
    \draw[black, ultra thick] (0,0) -- (0,-0.8);
    \draw[black, ultra thick] (0,0.8) -- (0,-0.8);
\end{tikzpicture}}_{=D_2FD_2FF}\\
& \qquad \quad + \underbrace{\begin{tikzpicture}[baseline={0}]
    \filldraw[black] (-0.4,0.4) circle (2pt);
    \filldraw[black] (0,-0.4) circle (2pt);
    \filldraw[black] (0.4,0.4) circle (2pt);
    \draw[black, dashed, ultra thick] (-0.4,0.4) -- (0,-0.4);
    \draw[black, dashed, ultra thick] (0.4,0.4) -- (0,-0.4);
\end{tikzpicture}}_{=D_{11}F[F,F]} + \underbrace{\begin{tikzpicture}[baseline={0}]
    \filldraw[black] (-0.4,0.4) circle (2pt);
    \filldraw[black] (0,-0.4) circle (2pt);
    \filldraw[black] (0.4,0.4) circle (2pt);
    \draw[black, dashed, ultra thick] (-0.4,0.4) -- (0,-0.4);
    \draw[black, ultra thick] (0.4,0.4) -- (0,-0.4);
\end{tikzpicture}}_{=D_{21}F[F,F]} + \underbrace{\begin{tikzpicture}[baseline={0}]
    \filldraw[black] (-0.4,0.4) circle (2pt);
    \filldraw[black] (0,-0.4) circle (2pt);
    \filldraw[black] (0.4,0.4) circle (2pt);
    \draw[black, ultra thick] (-0.4,0.4) -- (0,-0.4);
    \draw[black, dashed, ultra thick] (0.4,0.4) -- (0,-0.4);
\end{tikzpicture}}_{=D_{12}F[F,F]} + \underbrace{\begin{tikzpicture}[baseline={0}]
    \filldraw[black] (-0.4,0.4) circle (2pt);
    \filldraw[black] (0,-0.4) circle (2pt);
    \filldraw[black] (0.4,0.4) circle (2pt);
    \draw[black, ultra thick] (-0.4,0.4) -- (0,-0.4);
    \draw[black, ultra thick] (0.4,0.4) -- (0,-0.4);
\end{tikzpicture}}_{=D_{22}F[F,F]}\ .
\end{align*}

Given such an edge-colored tree $\tau$, we denote by $\mathcal{F}(\tau)(y)$ the corresponding elementary differential obtained by differentiation and insertion of $F$ as described above, e.g.,
\begin{align*}
    \mathcal{F} \left(\begin{tikzpicture}
    \filldraw[black] (0,0) circle (2pt);
\end{tikzpicture}\right)(y) &= F(y,y),\\
    \mathcal{F} \left(\begin{tikzpicture}[baseline={0}]
    \filldraw[black] (0,0.4) circle (2pt);
    \filldraw[black] (0,-0.4) circle (2pt);
    \draw[black, dashed, ultra thick] (0,0.4) -- (0,-0.4);
\end{tikzpicture}\right)(y) &= D_1F(y,y) F(y,y).
\end{align*}
Observe that, by symmetry of partial differentiation, the trees 
$$ \underbrace{\begin{tikzpicture}[baseline={0}]
    \filldraw[black] (-0.4,0.4) circle (2pt);
    \filldraw[black] (0,-0.4) circle (2pt);
    \filldraw[black] (0.4,0.4) circle (2pt);
    \draw[black, dashed, ultra thick] (-0.4,0.4) -- (0,-0.4);
    \draw[black, ultra thick] (0.4,0.4) -- (0,-0.4);
\end{tikzpicture}}_{=D_{21}F[F,F]} ,\ \underbrace{\begin{tikzpicture}[baseline={0}]
    \filldraw[black] (-0.4,0.4) circle (2pt);
    \filldraw[black] (0,-0.4) circle (2pt);
    \filldraw[black] (0.4,0.4) circle (2pt);
    \draw[black, ultra thick] (-0.4,0.4) -- (0,-0.4);
    \draw[black, dashed, ultra thick] (0.4,0.4) -- (0,-0.4);
\end{tikzpicture}}_{=D_{12}F[F,F]} $$
have the same elementary differential. We also have the usual internal symmetries of rooted trees. We define such trees to be equivalent and introduce a symmetry factor $\alpha(\tau)$ counting the elements of the equivalence class of such trees, e.g.,
$$ \alpha \left(\begin{tikzpicture}[baseline={0}]
    \filldraw[black] (-0.4,0.4) circle (2pt);
    \filldraw[black] (0,-0.4) circle (2pt);
    \filldraw[black] (0.4,0.4) circle (2pt);
    \draw[black, dashed, ultra thick] (-0.4,0.4) -- (0,-0.4);
    \draw[black, ultra thick] (0.4,0.4) -- (0,-0.4);
\end{tikzpicture} \right) = 2. $$

To be more precise, we introduce a recursive definition of edge-colored rooted trees, defining such a tree in terms of its children. Let $\tau$ be an edge-colored rooted tree and let the order $|\tau|$ be the number of its nodes. Let $\tau_1,\dots,\tau_m$ be the trees obtained by removing the root node and the edges connecting these trees to the root node. Let the edges connecting the root node of $\tau$ to $\tau_1,\dots,\tau_m$ have colors $a_1,\dots,a_m \in \{1,2\}$, respectively. Then, we express the tree $\tau$ as
$$ \tau = [\tau_1|_{a_1},\dots,\tau_m|_{a_m}]. $$
The elementary differential of $\tau$ can then be defined recursively as
\begin{subequations}\label{eq:elem-diff-recursive}
\begin{align}
    \mathcal{F} \left(\begin{tikzpicture}
    \filldraw[black] (0,0) circle (2pt);
\end{tikzpicture}\right)(y) &= F(y,y),\\
	\mathcal{F}(\tau)(y) &= D_{a_1\dots a_m}F(y,y) [\mathcal{F}(\tau_1)(y),\dots,\mathcal{F}(\tau_m)(y)] \text{ for } \tau=[\tau_1|_{a_1},\dots,\tau_m|_{a_m}],
\end{align} 
\end{subequations}
where $D_{a_1\dots a_m}F(y,y): X \times \dots \times X \rightarrow X$ denotes the $m-$multilinear mapping given by the $m^{th}$ derivative of $F(y,y)$ with respect to its $a_1,\dots,a_m$ arguments (for a discussion of higher order derivatives as multilinear maps, see \cite{cartan1971}).

The symmetry factor $\alpha(\tau)$ can be understood from the symmetry of the multilinear operator $D_{a_1\dots a_m}F$. In determining order conditions, we only want to consider independent elementary differentials, so we only wish to sum over distinct trees corresponding to different elementary differentials. The symmetry factor can be computed recursively as $\alpha \left(\begin{tikzpicture}
    \filldraw[black] (0,0) circle (2pt);
\end{tikzpicture}\right) = 1$ and 
\begin{align*}
\alpha (\tau) &= \begin{pmatrix} |\tau|-1 \\ |\tau_1|,\dots,|\tau_m| \end{pmatrix} \alpha(\tau_1)\cdots \alpha(\tau_m) \frac{1}{\mu_1!\cdots \mu_{j_m}!},
\end{align*}
where $\tau = [\tau_1|_{a_1},\dots,\tau_m|_{a_m}]$ and $\mu_1,\dots,\mu_{j_m}$ count the number of mutually equal trees with the same root edge-coloring $\tau_1|_{a_1},\dots,\tau_m|_{a_m}$, whose representatives we index as $\{1,\dots,j_m\}$. To see this, we consider labeling each non-root node with a number $1,\dots,|\tau|-1$; the multinomial coefficient gives the number of possible partitions of $1,\dots,|\tau|-1$ to the trees $\tau_1,\dots,\tau_m$. For each tree $\tau_i$, there are $\alpha(\tau_i)$ ways of assigning these labels. Finally, we divide by $\mu_1!\cdots \mu_{j_m}!$ since permutations of equal trees with the same root edge-coloring $\tau_1|_{a_1},\dots,\tau_m|_{a_m}$ do not change the labeling. This is derived in \cite{araujo1997}, with the slight modification that, in our case, a recursively defined edge-colored rooted tree $\tau = [\tau_1|_{a_1},\dots,\tau_m|_{a_m}]$ keeps track of the root edge labelings, so we appropriately modify the definition of the counting factors $\mu_1,\dots,\mu_{j_m}$. Throughout, when we sum over trees of order $q$, we understand that we are summing over equivalence classes of trees of order $q$, denoted $|\tau|=q$. In other words, a sum over non-identified trees of order $q$ can be expressed as a sum over equivalence classes of trees of order $q$ by introducing the symmetry factor $\alpha$, i.e.,
$$ \sum_{\text{non-identified trees of order } q} (\dots) = \sum_{|\tau|=q} \alpha(\tau) (\dots), $$
given that the summand $(\dots)$ is invariant on equivalence classes. In particular, the elementary differentials are invariant on equivalence classes, due to the symmetry of the multilinear operator $D_{a_1\dots a_m}F$.

Then, we have that the $q^{th}$ derivative can be expressed as a sum over edge-colored trees of order $q$.

\begin{proposition}
    For the exact solution $y(t)$, its $q^{th}$ derivative is given by
    \begin{equation}\label{eq:exact-sol-derivative-tree}
        \frac{d^q}{dt^q}y\Big|_{t_0} = \sum_{|\tau|=q} \alpha(\tau) \mathcal{F}(\tau)(y_0).
    \end{equation}
\end{proposition}
Note that the stages can be viewed of as functions of $h$, $Y_i = Y_i(h)$, implicitly defined from \eqref{eq:nprk-general}. Now, we define the following function of $h$, $g_{ij}(h) := h F(Y_i(h),Y_j(h))$. By the Leibniz rule (see Section III.1.1, Equation (1.8), of \cite{Hairer.2002}), we have
\begin{equation}\label{eq:q-derivative-g}
    \frac{d^q}{dh^q} g_{ij}\Big|_{h=0} = q \frac{d^{q-1}}{dh^{q-1}} F(Y_i,Y_j) \Big|_{h=0}.
\end{equation}
The factor of $q$ appearing on the right hand side of \eqref{eq:q-derivative-g} will produce additional integer factors in the Taylor expansion for the numerical solution about $h=0$. To account for these additional integer factors, noting that $q = |\tau|$, we define the density of a tree $\gamma(\tau)$ as the product of the order of $\tau$ with all orders of trees that appear if roots are successively removed, i.e.,
\begin{align*}
\gamma\left(\begin{tikzpicture}
    \filldraw[black] (0,0) circle (2pt);
\end{tikzpicture}\right) &= 1,\\
\gamma(\tau) &= |\tau| \gamma(\tau_1)\cdots \gamma(\tau_m),
\end{align*}
where $\tau = [\tau_1|_{a_1},\dots,\tau_m|_{a_m}]$.

Furthermore, the above formula \eqref{eq:q-derivative-g} allows us to inductively construct the elementary weights associated to a tree. To construct the elementary weights, we first note that \eqref{eq:q-derivative-g} provides an expression for the derivative of the internal stages of an NPRK method \eqref{eq:nprk-general}:
$$ \frac{d^q}{dh^q} Y_i = q\sum_{kl} a_{ikl} \frac{d^{q-1}}{dh^{q-1}} F(Y_k,Y_l)\Big|_{h=0}, $$
which is seen by taking the $q^{th}$ derivative of \eqref{eq:nprk-general} at $h=0$ and substituting \eqref{eq:q-derivative-g}. We associate to $F(Y_k,Y_l)$ a multi-indexed node $\, \begin{tikzpicture} [baseline={-2}] \filldraw[black] (0,0) circle (2pt) node[anchor=west]{$kl$};
\end{tikzpicture}.$ Differentiation of $F(Y_k,Y_l)$ with respect to the first argument produces a factor  $\sum_{uv}a_{kuv}(...)$ (here $...$ denotes differentials evaluated at $(Y_u,Y_v)$) and similarly $\sum_{uv}a_{luv}(...)$ for the second argument. Letting a tree $\tau$ have root labeled $ij$, we define the \textit{indexed elementary weight} $\phi_{ij}(\tau)$ as follows: label the remaining nodes in the tree with a multi-index; between a pair of nodes labeled $kl$ below and $uv$ above, if they are connected by a dashed edge, we obtain a factor $a_{kuv}$ and similarly, $a_{luv}$ for a solid edge; once we have traversed through all non-rooted nodes in the tree, we sum over the non-rooted indices. That is, differentiation with respect to the first argument produces a factor $a_{kuv}$ where the first index $k$ corresponds to the first component of the multi-index of the lower node and differentiation with respect to the second argument produces a factor $a_{luv}$ where the first index $l$ corresponds to the second component of the multi-index of the lower node. For example,
\begin{align*}
    \phi_{ij}\left(\begin{tikzpicture}[baseline={0}]
    \filldraw[black] (-0.4,0.4) circle (2pt) node[anchor=east]{$kl$};
    \filldraw[black] (0,-0.4) circle (2pt) node[anchor=west]{$ij$};
    \filldraw[black] (0.4,0.4) circle (2pt) node[anchor=west]{$uv$};
    \draw[black, dashed, ultra thick] (-0.4,0.4) -- (0,-0.4);
    \draw[black, ultra thick] (0.4,0.4) -- (0,-0.4);
\end{tikzpicture}\right) &= \sum_{kluv} a_{ikl}a_{juv}, 
\hspace{8ex}
    \phi_{ij}\left(\begin{tikzpicture}[baseline={0}]
    \filldraw[black] (0,0.8) circle (2pt) node[anchor=west]{$uv$};
    \filldraw[black] (0,0) circle (2pt) node[anchor=west]{$kl$};
    \filldraw[black] (0,-0.8) circle (2pt) node[anchor=west]{$ij$};
    \draw[black, ultra thick] (0,0) -- (0,-0.8);
    \draw[black, dashed, ultra thick] (0,0) -- (0,0.8);
\end{tikzpicture}\right) = \sum_{kluv}a_{jkl}a_{kuv}
\end{align*}
(note: the multi-indices $kl$ and $uv$ in the arguments of $\phi_{ij}$ above are dummy indices only shown for conceptual clarity). The preceding discussion produces a formula for the $q^{th}$ derivative of $g_{ij}$ at $h=0$,
$$ \frac{d^q}{dh^q}  g_{ij}\Big|_{h=0} = \sum_{|\tau|=q} \alpha(\tau) \gamma(\tau) \phi_{ij}(\tau) \mathcal{F}(\tau)(y_0). $$
This yields a formula for the $q^{th}$ derivative of $y_1 = y_0 + \sum_{ij}b_{ij}g_{ij}$ at $h=0$,
\begin{equation}\label{eq:num-sol-derivative-tree}
    \frac{d^q}{dh^q} y\Big|_{h=0} = \sum_{|\tau|=q} \alpha(\tau) \gamma(\tau) \sum_{ij}b_{ij}\phi_{ij}(\tau) \mathcal{F}(\tau)(y_0).
\end{equation}
We define the \textit{elementary weight} of a tree as $\Phi(\tau) = \sum_{ij}b_{ij}\phi_{ij}(\tau)$. Comparing equations \eqref{eq:exact-sol-derivative-tree} and \eqref{eq:num-sol-derivative-tree}, we obtain the order conditions for an NPRK method:
\begin{theorem}[NPRK Order Conditions]\label{thm:2-nprk-order-cond}
    An NPRK method has order $p$ if and only if
    \begin{equation}\label{eq:2-nprk-order-cond}
        \Phi(\tau) = \frac{1}{\gamma(\tau)} 
    \end{equation}
    for all trees $\tau$ such that $|\tau|\leq p.$
\begin{proof}
    From the above discussion, we have already proven the ``if" direction, so it remains to prove the ``only if" direction. As is standard (e.g., see Section II.2, Theorem 2.13, of \cite{Hairer.1993}), it suffices to show that the elementary differentials are independent. To show that they are independent, it suffices to show that for every  tree $\tau$ (more precisely, for the equivalence class corresponding to $\tau$), there exists a partitioned ODE $y' = F(y,y)$ on $\mathbb{R}^k$, for some $k$, with initial condition $y(0)$ such that
    \begin{equation}\label{eq:tree-to-show}
        \mathcal{F}_k(\tau)(y(0)) = 1,
    \end{equation}
    where $\mathcal{F}_k(\tau)(y(0))$ denotes the $k^{th}$ component of $\mathcal{F}(\tau)(y(0))$, and that 
    \begin{equation}\label{eq:not-tree-to-show}
        \mathcal{F}_k(\sigma)(y(0)) = 0
    \end{equation}
    for all other trees $\sigma \neq \tau$. We adapt the proof of the analogous result for Runge--Kutta methods (see Section II.2, Theorem 2.13 and Exercise 4, of \cite{Hairer.1993}) to account for edge coloring.
    
    For a tree of order $q$, we define an ODE on $\mathbb{R}^q \ni y = (y_1,\dots,y_q)$ with zero initial condition $y(0) = (0,\dots,0)$ as follows. For the tree $\begin{tikzpicture} \filldraw[black] (0,0) circle (2pt); \end{tikzpicture}$ of order $1$, we define the ODE to be simply $y_1' = 1$. Let $\tau$ be a tree of order $q$, expressed in recursive form  
    $$ \tau=[\tau_{q-1}|_{a_1},\dots,\tau_{q-m}|_{a_m}],$$
    where of course the number of children $m$ is less than $q$. We define the ODE to be $y' = f(y)$ where the $q^{th}$ component of $f$ is
    \begin{equation}\label{eq:qth-component}
        f(y)_q = \prod_{k = 1}^m y_{q-k}.
    \end{equation}
    Let $A_1$ be the set of all indices $j \in \{1,\dots,m\}$ such that $a_j=1$ and similarly let $A_2$ be the set of all indices $j \in \{1,\dots,m\}$ such that $a_j=2$. Then, define the $q^{th}$ component of the partition, corresponding to \eqref{eq:qth-component}, to be
    \begin{equation}\label{eq:qth-partition}
        F_q(x,y) = \prod_{k \in A_1} x_{q-k} \prod_{k \in A_2} y_{q-k}.
    \end{equation}
    Now, recall that the elementary differential at $y(0)=0$ is given in recursive form as
    $$ \mathcal{F}(\tau)(0) = D_{a_1\dots a_m}F(0,0) [\mathcal{F}(\tau_{q-1})(0),\dots,\mathcal{F}(\tau_{q-m})(0)]. $$
    \sloppy Denoting the components of the multi-linear operator as $D_{a_1\dots a_m}F(0,0)$ as $[D_{a_1\dots a_m}F(0,0)]^i_{j_1, \dots, j_m}$, the only non-zero component is
    $$ [D_{a_1\dots a_m}F(0,0)]^q_{q-1, \dots, q-m} = 1 $$
    by construction. Now, we finish by induction. Assume that \eqref{eq:tree-to-show} holds for all $k < q$, then 
    \begin{align*}
        \mathcal{F}_q(\tau)(0) &= \sum_{j_1,\dots,j_m = 1}^q [D_{a_1\dots a_m}F(0,0)]^q_{j_1, \dots, j_m} [\mathcal{F}_{j_1}(\tau_{q-1})(0),\dots,\mathcal{F}_{j_m}(\tau_{q-m})(0)] \\
        &=  [D_{a_1\dots a_m}F(0,0)]^q_{q-1, \dots, q-m} [\mathcal{F}_{q-1}(\tau_{q-1})(0),\dots,\mathcal{F}_{q-m}(\tau_{q-m})(0)] = 1.
    \end{align*}
    A similar argument shows that $\mathcal{F}_q(\sigma)(0) = 0$ for all $\sigma \neq \tau$, which completes the proof.
\end{proof}
\end{theorem}

\begin{remark}
    Note that, in principle, there is an alternative method to derive the order conditions for NPRK methods by viewing them as PRK methods. Namely, for the ODE $y' = F(y,y), y(0) = y_0$, one can consider the equivalent partitioned ODE $y' = F(y,z)$, $z' = G(y,z)$, $y(0)= y_0 = z(0)$. Subsequently, one can use the order condition theory for PRK methods to derive the order conditions for NPRK methods.

    For example, we do this in \cite{nprk1} for a simplified class of sequentially-coupled NPRK methods with a highly sparse NPRK tableau. For an $s$ stage sequentially-coupled method, it can equivalently be expressed as an $s-1$ stage PRK method. However, for a fully dense NPRK tableau, the corresponding PRK method would require $s^2$ stages; furthermore, there are many constraints that arise due to the duplication of variables. As such, we find this approach unnatural and instead work directly in the NPRK formulation. Furthermore, as we will see in \Cref{section:relation-to-ark-order-cond}, this direct approach will naturally reveal a separation between additive and nonlinear order conditions. We use these nonlinear order conditions in \Cref{sec:example} to show how nonlinear coupling strength can be estimated. 
\end{remark}

\subsection{Generalization to $M$ partitions}\label{sec:order:n}
Analogous to how two-component ARK methods can be generalized to $M$ components \cite{cooper1980additive}, NPRK methods can be generalized to $M$ partitions; the order conditions follow in a conceptually similar manner to the $M=2$ case.

Let $M$ be a positive integer. Consider the initial value problem 
\begin{align*}
    y' = f(y), \quad y(t_0) = y_0,
\end{align*}
specified by a vector field $f: X \rightarrow X.$ We say a map
$$ F: \underbrace{X \times \dots \times X}_{M \text{ times}} \rightarrow X $$
is an $M-$partition of $f: X \rightarrow X$ if
$F(y,\dots,y) = f(y) \text{ for all } y \in X.$
Given an $M-$partition $F$ of $f$, the ansatz for an $s$-stage NPRK$_M$ method for the above initial value problem is
\begin{align}
    \begin{aligned}
        Y_{i_0} &= y_n + h\sum_{i_1\dots i_M=1}^s a_{i_0 i_1 \dots i_M} F(Y_{i_1},\dots,Y_{i_M}), \quad i_0=1,\ldots, s, \\
        y_{n+1} &= y_n + h\sum_{i_1\dots i_M=1}^s b_{i_1 \dots i_M} F(Y_{i_1},\dots,Y_{i_M}).
    \end{aligned}
    \label{eq:n-nprk-general}
\end{align}
The rank $M+1$ tensor $a_{i_0 i_1 \dots i_M}$ replaces the classical RK matrix $a_{ij}$ and the rank $M$ tensor $b_{i_1 \dots i_M}$ replaces the classical RK weight vector $b_i$. The choice $M=1$ reproduces classical RK methods and the choice $M=2$ reproduces NPRK methods as defined previously. 

Every NPRK$_M$ method has $M$ underlying RK schemes, where the $r^{th}$ scheme, $r=1,\dots,M$, is given by using a trivial partition in the $r^{th}$ argument, $F(Y_{i_1},\dots,Y_{i_M}) = f(Y_{i_r})$. The corresponding RK coefficients are
$$ a^{\{r\}}_{i_0i_r} = \sum_{i_1,\dots,\widehat{i_r},\dots,i_M=1}^s a_{i_0i_1\dots i_M}, \quad b^{\{r\}}_{i_r} =  \sum_{i_1,\dots,\widehat{i_r},\dots,i_M=1}^s b_{i_1\dots i_M} $$
(where $\widehat{i_r}$ denotes omission of that index in the sum).

To derive the order conditions, we note that the time derivative is given by
$$ \frac{d}{dt}(\ ) = D_1(\ )\cdot F + \ldots + D_M(\ )\cdot F, $$
where again $D_i$ denotes differentiation with respect to the $k^{th}$ argument, $k=1,\dots,M$; for example,
$$ y'' = \frac{d}{dt}F(y,\dots,y) = D_1 F(y,\dots,y) \cdot F(y,\dots,y) + \ldots + D_M F(y,\dots,y) \cdot F(y,\dots,y). $$

The order conditions can be obtained analogously to the $M=2$ case; namely, we consider edge-colored rooted trees where each edge can be colored by one of $M$ colors $1,\dots,M$; we represent this graphically by writing the color adjacent to the edge. Differentiation is then given by grafting onto a rooted tree a colored edge connected to a node:
$$ D_1(\ )\cdot F = \begin{tikzpicture}[baseline={0}]
    \filldraw[black] (0,0.4) circle (2pt);
    \draw[black, ultra thick] (0,0.4) -- (0,-0.4) node[anchor=west,midway]{$1$};
\end{tikzpicture}\ ,\dots,\ D_M(\ )\cdot F = \begin{tikzpicture}[baseline={0}]
    \filldraw[black] (0,0.4) circle (2pt);
    \draw[black, ultra thick] (0,0.4) -- (0,-0.4) node[anchor=west,midway]{$M$};
\end{tikzpicture}. $$
For example, for two colors $a,b \in \{1,\dots,M\}$,
\begin{align*}
    \begin{tikzpicture}[baseline={0}]
    \filldraw[black] (-0.4,0.4) circle (2pt);
    \filldraw[black] (0,-0.4) circle (2pt);
    \filldraw[black] (0.4,0.4) circle (2pt);
    \draw[black, ultra thick] (-0.4,0.4) -- (0,-0.4) node[anchor=east,midway]{$a$};
    \draw[black, ultra thick] (0.4,0.4) -- (0,-0.4) node[anchor=west,midway]{$b$};
\end{tikzpicture} &= D_{ab}F[F,F],
\hspace{8ex}
    \begin{tikzpicture}[baseline={0}]
    \filldraw[black] (0,0.8) circle (2pt) ;
    \filldraw[black] (0,0) circle (2pt) ;
    \filldraw[black] (0,-0.8) circle (2pt);
    \draw[black, ultra thick] (0,0) -- (0,-0.8) node[anchor=west,midway]{$a$};
    \draw[black, ultra thick] (0,0) -- (0,0.8) node[anchor=west,midway]{$b$};
\end{tikzpicture} = D_aF D_bF F.
\end{align*}

We define the order $|\tau|$, density $\gamma(\tau)$, symmetry factor $\alpha(\tau)$, and elementary differential $\mathcal{F}(\tau)(y)$ of a tree analogous to Section \ref{subsection: tree order conditions}. That is, the definitions are formally the same, with the modification that a recursively-defined edge-colored rooted tree $\tau = [\tau_1|_{a_1},\dots,\tau_m|_{a_m}]$ now has colorings $a_1,\dots,a_m$ valued in $\{1,\dots,M\}$.

The indexed elementary weight of a tree, $\phi_{i_1\dots i_M}(\tau)$, is defined analogously: with the root of $\tau$ labeled with the $M$-multi-index $i_1\dots i_M$, label every other node of $\tau$ with an $M$-multi-index; given two nodes with the lower node indexed by $j_1\dots j_M$ and the upper node indexed by $k_1\dots k_M$ connected by an edge of color $a$, the elementary weight receives a factor
$$ a_{j_a k_1 \dots k_M}(\dots), $$
i.e., the index $j_a$ is the $a^{th}$ index of the lower node, corresponding to the color $a$ of the edge connecting the two nodes. For example,
\begin{align*}
    \phi_{i_1\dots i_M}\left(\begin{tikzpicture}[baseline={0}]
    \filldraw[black] (-0.8,0.8) circle (2pt) node[anchor=east]{$j_1\dots  j_M$};
    \filldraw[black] (0,-0.8) circle (2pt) node[anchor=west]{$i_1\dots  i_M$};
    \filldraw[black] (0.8,0.8) circle (2pt) node[anchor=west]{$k_1\dots  k_M$};
    \draw[black, ultra thick] (-0.8,0.8) -- (0,-0.8) node[anchor=east,midway]{$a$};
    \draw[black, ultra thick] (0.8,0.8) -- (0,-0.8) node[anchor=west,midway]{$b$};
\end{tikzpicture}\right) &= \sum_{j_1\dots j_M} \sum_{k_1\dots k_M}  a_{i_a j_1\dots j_M} a_{i_b k_1\dots k_M}, \\
    \phi_{i_1\dots i_M}\left(\begin{tikzpicture}[baseline={0}]
    \filldraw[black] (0,0.8) circle (2pt) node[anchor=west]{$k_1\dots k_M$};
    \filldraw[black] (0,0) circle (2pt) node[anchor=west]{$j_1\dots j_M$};
    \filldraw[black] (0,-0.8) circle (2pt) node[anchor=west]{$i_1\dots i_M$};
    \draw[black, ultra thick] (0,0) -- (0,-0.8) node[anchor=east,midway]{$a$};
    \draw[black, ultra thick] (0,0) -- (0,0.8) node[anchor=east,midway]{$b$};
\end{tikzpicture}\right) &= \sum_{j_1\dots j_M} \sum_{k_1\dots k_M} a_{i_a j_1 \dots j_M}  a_{j_b k_1 \dots k_M}.
\end{align*}
Defining the \textit{elementary weight}
$$ \Phi(\tau) = \sum_{i_1\dots i_M}b_{i_1\dots i_M}\phi_{i_1\dots i_M}(\tau), $$
an analogous argument to Section \ref{subsection: tree order conditions} gives the order conditions for an NPRK$_M$ method: 
\begin{theorem}[NPRK$_M$ Order Conditions]\label{thm:n-nprk-order-cond}
    An NPRK$_M$ method has order $p$ if and only if
    \begin{equation}\label{eq:n-nprk-order-cond}
        \Phi(\tau) = \frac{1}{\gamma(\tau)}
    \end{equation}
    for all trees $\tau$ such that $|\tau|\leq p$.
\end{theorem}

\subsection{Enumeration and computation of NPRK$_M$ order conditions}\label{sec:order:enum}
The enumeration of the NPRK$_M$ order conditions is given by counting the number of $M-$edge-colored rooted trees at each order. For $M$ up to $10$, the enumeration of such trees is given in \cite{Foissy.2021}, where such trees are referred to as $M-$typed $1-$decorated rooted trees (the type refers to coloring of edges and the decoration refers to coloring of nodes). 

Let $\sigma^{[M]}_i$ denote the number of NPRK$_M$ order conditions of order $i$. Following \cite{Kennedy.2003tv4}, let $M\alpha^{[M]}_i$ denote the number of $ARK_M$ order conditions of order $i$, where $\alpha^{[M]}_i$ is defined from the generating function
$$ \sum_{i=1}^\infty \alpha^{[M]}_i x^{i-1} = \prod_{i=1}^\infty (1-x^i)^{-M\alpha^{[M]}_i}$$
(here, $x$ is an indeterminate in a formal power series).

We will additionally count the number of coupling conditions for an NPRK$_M$ method, i.e., order conditions that do not correspond to order conditions for the $M$ underlying RK methods; we discuss this further in \Cref{section:relation-to-ark-order-cond}. This is straightforward as the order conditions for the $M$ underlying RK methods arise from only considering trees whose edges have a single color. Thus, for an NPRK$_M$ method, there are $M\alpha^{[1]}_i$ order conditions for the underlying methods. Denote the number of coupling conditions for an NPRK$_M$ method by $\widetilde{\sigma}^{[M]}_i$ ($M \geq 2$). Then,
$$ \widetilde{\sigma}^{[M]}_i = \sigma^{[M]}_i - M\alpha^{[1]}_i.$$
The enumeration of both the total number of order conditions and the number of coupling conditions is shown in Table \ref{table:order-cond-enumeration} for the number of partitions $M$ up to $5$ and order $i$ up to $8$.
\begin{table}[H]
\caption{Enumeration of NPRK$_M$ order conditions for up to 5 partitions and order up to $8$.}\label{table:order-cond-enumeration}
\begin{center}
\begin{tabular}{|c |c |c c c c c c c c|} 
 \hline
 $M$ & order $i:$ & 1 & 2 & 3 & 4 & 5 & 6 & 7 & 8 \\
 \hline
 \hline
 $1$ & $\sigma^{[1]}_i$ & 1 & 1 & 2 & 4 & 9 & 20 & 48 & 115 \\
 \hline 
 $2$ & $\sigma^{[2]}_i$ & 1 & 2 & 7 & 26 & 107 & 458 & 2058 & 9498 \\
     & $\widetilde{\sigma}^{[2]}_i$ & 0 & 0 & 3 & 18 & 89 & 418 & 1962 & 9268 \\ \hline
 $3$ & $\sigma^{[3]}_i$ & 1 & 3 & 15 & 82 & 495 & 3144 & 20875 & 142773 \\
     & $\widetilde{\sigma}^{[3]}_i$ & 0 & 0 & 9 & 70 & 468 & 3084 & 20731 & 142428 \\ \hline     
 $4$ & $\sigma^{[4]}_i$ & 1 & 4 & 26 & 188 & 1499 & 12628 & 111064 & 1006840 \\
     & $\widetilde{\sigma}^{[4]}_i$ & 0 & 0 & 18 & 172 & 1463 & 12548 & 110872 & 1006380 \\ \hline   
 $5$ & $\sigma^{[5]}_i$ & 1 & 5 & 40 & 360 & 3570 & 37476 & 410490 & 4635330 \\
     & $\widetilde{\sigma}^{[5]}_i$ & 0 & 0 & 30 & 340 & 3525 & 37376 & 410250 & 4634755  \\
     \hline
\end{tabular}
\end{center}
\end{table}

\begin{remark}\label{remark:nprk-ark-tree-correspond}
Observe that the number of NPRK$_M$ order conditions is $1/M$ times the number of ARK$_M$ order conditions, i.e.,
$$ \sigma^{[M]}_i = \alpha^{[M]}_i. $$
This can be explained by the fact that each NPRK$_M$ tree corresponds to $M$ ARK$_M$ trees. To see this, given an $M-$edge-colored tree, for the root node, create $M$ trees by coloring the root with one of $M$ colors; for every other node, color them with the color of the edge below it and finally remove all of the edge colors (note that this procedure is purely for enumeration; it does not map NPRK elementary differentials correctly to ARK elementary differentials). By this procedure, the set of all equivalence classes of $M-$NPRK trees corresponds to the set of all equivalence classes of $\text{ARK}_M$ trees with a fixed root color. In the reverse direction, fixing the root color of an $\text{ARK}_M$ tree, given an $\text{ARK}_M$ tree with that root color, remove the coloring of the root node, color each edge by the color of the node above the edge, and finally remove all node colorings, producing an $M-$NPRK tree.

In terms of tableaux coefficients, the reduced number of order conditions for NPRK vs ARK methods, despite the larger number of mixed differentials to account for, can be seen as arising from the tensorial nature of the NPRK tableaux, which we discuss further in \Cref{section:relation-to-ark-order-cond}.
\end{remark}

\textbf{Computation of Order Conditions. }We will now present an algorithm for computing the order conditions. To be concrete, we will consider the case of $M=2$ partitions, although a higher number of partitions follows similarly. We know each NPRK order condition can be obtained by an equivalence class of NPRK trees, as discussed in \Cref{subsection: tree order conditions}. To choose a representative of each equivalence class, we begin with a canonical representative of each equivalence class of $\text{ARK}_2$ trees given by lexicographic ordering of the level and color sequences, as described in \cite{ketcheson.2023}. By considerng ARK$_2$ trees with a fixed root coloring, the procedure described in Remark \ref{remark:nprk-ark-tree-correspond} puts these trees in correspondence with a canonical representatives of NPRK$_2$ trees. We utilize the Julia package \verb|RootedTrees.jl| \cite{ranocha2019rootedtrees} to generate the canonical representative of each tree $\tau$, whose canonical representative is specified by a level sequence $L(\tau)$ and a color sequence $C(\tau)$ describing the level of each node in lexicographic order and the color of each node, respectively. For the $k^{th}$ element of the level sequence, $L(\tau)_k$, we denote by $\text{Parent}(\tau,k)$ the index of the parent node of node $k$, i.e., the largest index $p$ in $L(\tau)$ such that $L(\tau)_k - 1 = L(\tau)_p$.

\begin{remark}
Note that \verb|RootedTrees.jl| colors canonical bicolored rooted trees with $0$ and $1$; we will let the colors $0$ and $1$ correspond to differentiation with respect to the first and second arguments, respectively. Furthermore, since we only consider ARK$_2$ trees with a fixed root color, we will just consider canonical trees where the root is colored by $0$, as a matter of convention. As we will see in \Cref{alg:two-nprk-order} below, the first index of the color sequence is never used. 
\end{remark}

We then apply the method of obtaining the order conditions described in \Cref{subsection: tree order conditions}. Given a canonically represented tree as described above, we label each node with indices $i_kj_k$ where $k$ is the integer where the node appears in the level sequence. For an upper node with indices $i_kj_k$ and a lower node with indices $i_lj_l$ connected by an edge of color $b \in \{0,1\}$, we obtain a factor $a_{i_l i_kj_k}$ if $b = 0$ or $a_{j_l i_kj_k}$ if $b=1$. We thus loop through all non-rooted nodes in the tree, multiply the factors obtained from them, multiply by $b_{i_1j_1}$ where $i_1j_1$ are the indices of the root node, and sum over all indices to obtain the left-hand-side of the order condition, i.e., the elementary weight. The right hand side of the order condition is given by the reciprocal of the density of the tree, $1/\gamma(\tau)$, which can also be obtained from \verb|RootedTrees.jl|. This algorithm is summarized in \Cref{alg:two-nprk-order}.

\begin{algorithm}[H]
\caption{Order Condition for a $2-$NPRK Tree}\label{alg:two-nprk-order}
\begin{algorithmic}
\Require $a_{ijk}, b_{ij}$, $i,j,k=1,\dots,s$
\Require Canonical tree $\tau$, specified by (lexicographically ordered) level sequence $L(\tau)$, color sequence $C(\tau)$, density $\gamma(\tau)$.
\State $N \gets \text{Length}(L(\tau))$
\State $\text{sum} \gets 0$
\For{$i_1,j_1,\dots,i_N,i_N = 1$; $i_1,j_1,\dots,i_N,j_N \leq s$; $i_1,j_1,\dots,i_N,j_N\text{++}$}
\State $\text{prod} \gets 1$
\For{$k = N$; $k \geq 2$; $k\text{-}\text{-}$}
     \If{$C(\tau)_k == 0$}
          \State $\text{prod} \gets \text{prod}*a_{i_{\text{Parent}(\tau,k)}i_kj_k}$
     \Else
          \State $\text{prod} \gets \text{prod}*a_{j_{\text{Parent}(\tau,k)}i_kj_k}$
     \EndIf
\EndFor
\State $\text{sum} \gets \text{sum} + b_{i_1j_1}*\text{prod}$
\EndFor
\State return $\text{sum} == 1/\gamma(\tau)$
\end{algorithmic}
\end{algorithm}

In addition to a numerical implementation of \Cref{alg:two-nprk-order} in Julia, we also implemented a symbolic version of this algorithm in \verb|Mathematica| to generate the order conditions symbolically. For example, the third-order and fourth-order conditions for $M=2$ are shown in \Cref{app:order-list}, Equations \eqref{eq:n2order3list} and \eqref{eq:n2order4list}, respectively. Third-order conditions for $M=3$ are shown in \eqref{eq:n3order3list}. Both the numerical and symbolic implementations are available at \cite{githubNPRK2025}.

\section{Relation to additive order conditions}\label{section:relation-to-ark-order-cond}

As we have seen, an additive partition for an NPRK method results in an ARK method. In this section, we will relate the previously obtained order conditions for NPRK methods to order conditions for ARK methods (see, for example, \cite{araujo1997, Kennedy.2003tv4} and Section III.2 of \cite{Hairer.2002}). Note this discussion equally applies to PRK methods, as they can be expressed as ARK methods.

Note that we have shown that an NPRK method with tableaux $(a_{ijk}, b_{ij})$ with an additive partition reduces to a two-component ARK method, with the pair of tableau $(a^{\{1\}},b^{\{1\}},c^{\{1\}})$ and $(a^{\{2\}},b^{\{2\}},c^{\{2\}})$ (see \Cref{sec:underlying-rk-ark}). We will now state a partial converse.
\begin{proposition}
Consider a two-component ARK method with the pair of tableaux  $(a^{\{1\}},b^{\{1\}},c^{\{1\}})$ and $(a^{\{2\}},b^{\{2\}},c^{\{2\}})$. Additionally, assume that both tableaux are at least first-order $\sum_ib^{\{1\}}_i=1=\sum_jb^{\{2\}}_j$, have the same number of stages $s$, and satisfy $c^{\{1\}}=c^{\{2\}} = c$. Then, the method can be expressed as an NPRK method with an additive  partition with coefficients
\begin{subequations}
\begin{align}
a_{ijk} &= \frac{a^{\{1\}}_{ij}}{s} + \frac{a^{\{2\}}_{ik}}{s} - \frac{c_i}{s^2}, \label{eq:ark-to-nprk-coeff-a} \\
b_{ij} &= \frac{b^{\{1\}}_{i}}{s} + \frac{b^{\{2\}}_{j}}{s} - \frac{1}{s^2}. \label{eq:ark-to-nprk-coeff-b}
\end{align}
\end{subequations}
\begin{proof}
We simply have to check that the sums of $a$ (resp. $b$) over its second and third (resp. second) indices reproduces the pair of underlying RK tableaux.
\begin{align*}
\sum_{k=1}^s a_{ijk} &= \sum_{k=1}^s \left( \frac{a^{\{1\}}_{ij}}{s} + \frac{a^{\{2\}}_{ik}}{s} - \frac{c_i}{s^2} \right) = a^{\{1\}}_{ij} + \frac{c^{\{2\}}_i}{s} - \frac{c_i}{s} = a^{\{1\}}_{ij}, \\
\sum_{j=1}^s a_{ijk} &= \sum_{j=1}^s \left( \frac{a^{\{1\}}_{ij}}{s} + \frac{a^{\{2\}}_{ik}}{s} - \frac{c_i}{s^2} \right) = \frac{c^{\{1\}}_i}{s} + \frac{a^{\{2\}}_{ik}}{s} - \frac{c_i}{s} = a^{\{2\}}_{ik}, \\
\sum_j b_{ij} &= \sum_j \left(\frac{b^{\{1\}}_{i}}{s} + \frac{b^{\{2\}}_{j}}{s} - \frac{1}{s^2} \right) = b^{\{1\}}_i + \frac{1}{s} - \frac{1}{s} = b^{\{1\}}_i, \\
\sum_i b_{ij} &= \sum_i \left(\frac{b^{\{1\}}_{i}}{s} + \frac{b^{\{2\}}_{j}}{s} - \frac{1}{s^2} \right) = \frac{1}{s} + b^{\{2\}}_j - \frac{1}{s} = b^{\{2\}}_j.
\end{align*}
\end{proof}
\end{proposition}
\begin{remark}
Note that there can be different NPRK methods which have the same underlying ARK method. To see this, consider for example an ARK method where $b^{\{1\}}_i = b^{\{2\}}_i$, $i=1,\dots,s\geq 2$. Then, the ARK method can be expressed as an NPRK method with the above choice of $a$ and $b$, \eqref{eq:ark-to-nprk-coeff-a}-\eqref{eq:ark-to-nprk-coeff-b}. Note this choice of $b$ is not generally diagonal. On the other hand, it can also be expressed as an NPRK method with the above choice of $a$ and
\begin{equation}
b_{ij} = b^{\{1\}}_i \delta_{ij}, \label{eq:ark-to-nprk-coeff-diagonal}
\end{equation}
which is diagonal. These two possible choices correspond to different NPRK methods when applied to a general nonlinear partition but reduce to the same ARK method given an additive partition.
\end{remark}

In light of this discussion, we are able to relate ARK order conditions with NPRK order conditions with additive partitions. 

\textbf{Relation to ARK Order Conditions.} 
We will now investigate the relation of NPRK order conditions to ARK order conditions. In so doing, we will see that NPRK order conditions imply ARK order conditions but the converse is not true. In essence, this is because the linearly separable structure of ARK methods remove some elementary differentials appearing in a general NPRK method; namely, those containing $D_I F$, where $I=i_1\dots i_n$ is a multi-index such that at least two of its indices are different. In other words, nonlinear coupling of stages gives rise to new order conditions.

We will first consider the two-component case. Consider an ODE $y' = f(y) = f_1(y) + f_2(y)$ and define the partition
$$ F(y_1,y_2) = f_1(y_1) + f_2(y_2). $$

To understand how to relate the ARK and NPRK order conditions, let us start with the first-order condition. For an NPRK method, the first-order condition is
$$ \sum_{ij}b_{ij}=1. $$
On the other hand, for an ARK method with coefficients $(a^{\{1\}},b^{\{1\}},c)$ and $(a^{\{2\}},b^{\{2\}},c)$, there are two order conditions
$$ \sum_{i}b^{\{1\}}_{i} = 1 = \sum_{j}b^{\{2\}}_j. $$
The doubling of the number of order conditions arise from the fact that, for an NPRK method, $b^{\{1\}}$ and $b^{\{2\}}$ are not independent but rather, related via the tensor $b$. Namely,
$$ \sum_i b^{\{1\}}_i = \sum_{ij}b_{ij} = \sum_j b^{\{2\}}_j. $$
Thus, we only have to account for one first-order condition for an NPRK method, as the tensorial nature of the NPRK tableau automatically accounts for the other condition. We will visualize the process of converting an NPRK method to an ARK method by relating their rooted trees. For an NPRK method, we utilize edge-colored rooted trees as before; we refer to the dashed edge as having ``color 1" and the solid edge as having ``color 2". For an ARK method, it is standard to utilize node-colored rooted trees. We will color the nodes for an ARK method as
$$ \begin{tikzpicture}[baseline={0}] \draw[black] (0,0) circle (2pt);
\end{tikzpicture}\ ,\, \begin{tikzpicture}[baseline={0}] \filldraw[blue] (0,0) circle (2pt);
\end{tikzpicture}\,, $$
where the open node, referred to as ``color 1", corresponds to $f_1$ and the filled node, referred to as ``color 2", corresponds to $f_2$. We can thus visualize the rooted tree with one node for an NPRK method as decomposing into a pair of rooted trees for an ARK method,
$$\begin{tikzpicture}[baseline={0}] \filldraw[black] (0,0) circle (2pt) node[anchor=west]{$ij$};
\end{tikzpicture} \longrightarrow \begin{tikzpicture}[baseline={0}] \draw[black] (0,0) circle (2pt) node[anchor=west]{$i$};
\end{tikzpicture}\ ,\ \begin{tikzpicture}[baseline={0}] \filldraw[blue] (0,0) circle (2pt) node[anchor=west,black]{$j$};
\end{tikzpicture},$$
and thus, a decomposition of the NPRK order condition to two ARK order conditions.

To understand decomposing higher order NPRK trees into ARK trees, consider the derivative matrices of $F$ expressed in terms of $f_1$ and $f_2$
\begin{align*}
D_1F F &= D_1f_1 (f_1+f_2) = D_1f_1 f_1 + D_1f_1 f_2, \\
D_2F F &= D_2f_2 (f_1+f_2) = D_2f_2 f_1 + D_2f_2 f_2.
\end{align*}
This gives us a method of mapping an NPRK tree with two nodes connected by a colored edge to two ARK trees:
\begin{align*}
D_1FF &= \begin{tikzpicture}[baseline={0}]
    \filldraw[black] (0,0.4) circle (2pt);
    \filldraw[black] (0,-0.4) circle (2pt);
    \draw[black, dashed, ultra thick] (0,0.4) -- (0,-0.4);
\end{tikzpicture} \longrightarrow \begin{tikzpicture}[baseline={0}]
    \draw[black] (0,0.4) circle (2pt);
    \draw[black] (0,-0.4) circle (2pt);
    \draw[black] (0,0.36) -- (0,-0.36);
\end{tikzpicture} + \begin{tikzpicture}[baseline={0}]
    \filldraw[blue] (0,0.4) circle (2pt);
    \draw[black] (0,-0.4) circle (2pt);
    \draw[black] (0,0.36) -- (0,-0.36);
\end{tikzpicture} = D_1f_1 f_1 + D_1f_1 f_2,\\
D_2FF &= \begin{tikzpicture}[baseline={0}]
    \filldraw[black] (0,0.4) circle (2pt);
    \filldraw[black] (0,-0.4) circle (2pt);
    \draw[black, ultra thick] (0,0.4) -- (0,-0.4);
\end{tikzpicture} \longrightarrow \begin{tikzpicture}[baseline={0}]
    \draw[black] (0,0.4) circle (2pt);
    \filldraw[blue] (0,-0.4) circle (2pt);
    \draw[black] (0,0.36) -- (0,-0.36);
\end{tikzpicture} + \begin{tikzpicture}[baseline={0}]
    \filldraw[blue] (0,0.4) circle (2pt);
    \filldraw[blue] (0,-0.4) circle (2pt);
    \draw[black] (0,0.36) -- (0,-0.36);
\end{tikzpicture} = D_2 f_2f_1 + D_2f_2f_2.
\end{align*}
Namely, given an NPRK tree with two nodes connected by an edge of color $a \in \{1,2\}$, remove the coloring of the edge, color the lower node $a$ and form two trees that have upper nodes of color $1$ and $2$. Observe that the two second-order conditions for an NPRK method,
$$ \sum_{ijkl}b_{ij}a_{ikl} = \frac{1}{2} = \sum_{ijkl}b_{ij}a_{jkl}, $$
give the four second-order conditions for the underlying ARK method
$$ \sum_i b^{\{1\}}_i c^{\{2\}}_i = \sum_i b^{\{1\}}_i c^{\{1\}}_i = \frac{1}{2} = \sum_j b^{\{2\}}_j c^{\{2\}}_j =  \sum_j b^{\{2\}}_j c^{\{1\}}_j.$$
Analogous to the first-order case, the tensorial nature of the NPRK method means one only has to account for two second-order conditions. Of course, this is expected, if an NPRK method has order $p$, then the method with an additive partition must also have order $p$, since the order conditions derived for an NPRK method were independent of the partition. However, we continue to investigate the order conditions at higher order, as it will reveal new nonlinear order conditions that are not present in the ARK order conditions.

To generalize to the higher order case, we first introduce some terminology. We say that a node of a rooted tree is a \textit{leaf} if it does not connect to any nodes above it, i.e., has out-degree zero. We say that a node \textit{branches} if it connects to at least two nodes above it, i.e., has out-degree greater than or equal to two. We say that an edge-colored rooted tree is \textit{color-branching} if it contains a node that branches with at least two edges of different colors.

For the higher order case, we repeat the procedure above from top to bottom: we split an NPRK tree into multiple ARK trees, by starting with all possible colorings of the leaves. From there, the edges of the NPRK tree tell us how to color all of the lower nodes, as done above. At first, it may seem that this procedure is ill-defined, since a node that branches with two different colors would not have a well-defined color. However, such trees need not be considered, as the following proposition shows.

\begin{proposition}\label{prop:color-branching-differential}
The elementary differential of a color-branching tree, for an NPRK method with an additive partition, vanishes. 
\begin{proof}
Given a color-branching tree, one of the nodes branches with two edges of different colors. Thus, the elementary differential contains a factor
$$ \dots D_{12} F \dots \text{ or } \dots D_{21} F \dots $$
which vanishes, since $D_{12}F = 0 = D_{21}F$ for an additive partition $F(y_1,y_2) = f_1(y_1) + f_2(y_2)$. 
\end{proof}
\end{proposition}

This method is thus well-defined and allows us to see that all of the ARK order conditions of order $p$ are satisfied if the NPRK order conditions of order $p$ are satisfied (again, this is a trivial observation since if an NPRK method has order $p$, then the method with an additive partition must also have order $p$). For example, for third-order, we have six NPRK trees with nonzero differentials,
\begin{align*} \underbrace{\begin{tikzpicture}[baseline={0}]
    \filldraw[black] (0,0.8) circle (2pt);
    \filldraw[black] (0,0) circle (2pt);
    \filldraw[black] (0,-0.8) circle (2pt);
    \draw[black, dashed, ultra thick] (0,0) -- (0,-0.8);
    \draw[black, dashed, ultra thick] (0,0) -- (0,0.8);
\end{tikzpicture}}_{=D_1F D_1F F} &+ \underbrace{\begin{tikzpicture}[baseline={0}]
    \filldraw[black] (0,0.8) circle (2pt);
    \filldraw[black] (0,0) circle (2pt);
    \filldraw[black] (0,-0.8) circle (2pt);
    \draw[black, ultra thick] (0,0) -- (0,-0.8);
    \draw[black, dashed, ultra thick] (0,0) -- (0,0.8);
\end{tikzpicture}}_{=D_2FD_1FF} + \underbrace{\begin{tikzpicture}[baseline={0}]
    \filldraw[black] (0,0.8) circle (2pt);
    \filldraw[black] (0,0) circle (2pt);
    \filldraw[black] (0,-0.8) circle (2pt);
    \draw[black, dashed, ultra thick] (0,0) -- (0,-0.8);
    \draw[black, ultra thick] (0,0) -- (0,0.8);
\end{tikzpicture}}_{=D_1FD_2FF} + \underbrace{\begin{tikzpicture}[baseline={0}]
    \filldraw[black] (0,0.8) circle (2pt);
    \filldraw[black] (0,0) circle (2pt);
    \filldraw[black] (0,-0.8) circle (2pt);
    \draw[black, ultra thick] (0,0) -- (0,-0.8);
    \draw[black, ultra thick] (0,0.8) -- (0,-0.8);
\end{tikzpicture}}_{=D_2FD_2FF} \\ 
& + \underbrace{\begin{tikzpicture}[baseline={0}]
    \filldraw[black] (-0.4,0.4) circle (2pt);
    \filldraw[black] (0,-0.4) circle (2pt);
    \filldraw[black] (0.4,0.4) circle (2pt);
    \draw[black, dashed, ultra thick] (-0.4,0.4) -- (0,-0.4);
    \draw[black, dashed, ultra thick] (0.4,0.4) -- (0,-0.4);
\end{tikzpicture}}_{=D_{11}F[F,F]} +\ 2 \cancel{ \underbrace{\begin{tikzpicture}[baseline={0}]
    \filldraw[black] (-0.4,0.4) circle (2pt);
    \filldraw[black] (0,-0.4) circle (2pt);
    \filldraw[black] (0.4,0.4) circle (2pt);
    \draw[black, dashed, ultra thick] (-0.4,0.4) -- (0,-0.4);
    \draw[black, ultra thick] (0.4,0.4) -- (0,-0.4);
\end{tikzpicture}}_{=D_{21}F[F,F]}}  + \underbrace{\begin{tikzpicture}[baseline={0}]
    \filldraw[black] (-0.4,0.4) circle (2pt);
    \filldraw[black] (0,-0.4) circle (2pt);
    \filldraw[black] (0.4,0.4) circle (2pt);
    \draw[black, ultra thick] (-0.4,0.4) -- (0,-0.4);
    \draw[black, ultra thick] (0.4,0.4) -- (0,-0.4);
\end{tikzpicture}}_{=D_{22}F[F,F]}. 
\end{align*}
It is clear that the first four NPRK trees above correspond to eight ARK trees. The last two nonzero trees correspond to six ARK trees, 
\begin{align*} D_{11}F[F,F] &= \begin{tikzpicture}[baseline={0}]
    \filldraw[black] (-0.4,0.4) circle (2pt);
    \filldraw[black] (0,-0.4) circle (2pt);
    \filldraw[black] (0.4,0.4) circle (2pt);
    \draw[black, dashed, ultra thick] (-0.4,0.4) -- (0,-0.4);
    \draw[black, dashed, ultra thick] (0.4,0.4) -- (0,-0.4);
\end{tikzpicture} = \begin{tikzpicture}[baseline={0}]
    \draw[black] (-0.4,0.4) circle (2pt);
    \draw[black] (0,-0.4) circle (2pt);
    \draw[black] (0.4,0.4) circle (2pt);
    \draw[black] (-0.36,0.36) -- (0,-0.36);
    \draw[black] (0.36,0.36) -- (0,-0.36);
\end{tikzpicture} +\ 2\begin{tikzpicture}[baseline={0}]
    \filldraw[blue] (-0.4,0.4) circle (2pt);
    \draw[black] (0,-0.4) circle (2pt);
    \draw[black] (0.4,0.4) circle (2pt);
    \draw[black] (-0.36,0.36) -- (0,-0.36);
    \draw[black] (0.36,0.36) -- (0,-0.36);
\end{tikzpicture} + \begin{tikzpicture}[baseline={0}]
    \filldraw[blue] (-0.4,0.4) circle (2pt);
    \draw[black] (0,-0.4) circle (2pt);
    \filldraw[blue] (0.4,0.4) circle (2pt);
    \draw[black] (-0.36,0.36) -- (0,-0.36);
    \draw[black] (0.36,0.36) -- (0,-0.36);
\end{tikzpicture} \\ &= D^2f_1 [f_1,f_1] + 2D^2f_1[f_1,f_2] + D^2f_1[f_2,f_2] 
\end{align*}
and similarly $D_{22}F[F,F] = D^2f_2[f_1,f_1] + 2D^2f_2[f_2,f_1] + D^2f_2[f_2,f_2]$. Thus, the 6 non-vanishing third-order NPRK trees give rise to $14$ third-order ARK trees, which corresponds to the number of third-order conditions for two-component ARK methods \cite{Kennedy.2003tv4}.

The above construction shows that every NPRK condition corresponding to trees which are not color-branching can be expressed in terms of ARK order conditions. More precisely, Proposition \ref{prop:color-branching-differential} implies that any factor $a_{kuv}$ in the elementary weight of a non-color-branching tree will always have either $u$ or $v$ as a free index, i.e., there are no other factors in the elementary weight depending on either $u$ or $v$ since otherwise the tree must be color-branching. Once the final summation over all non-rooted indices is performed, this means $a_{kuv}$ can be expressed as either $a^{\{2\}}_{kv}$ or $a^{\{1\}}_{ku}$. Thus, all factors of $a$ appearing in the order condition of a non-color-branching tree can be expressed as either $a^{\{1\}}$ or $a^{\{2\}}$. Finally, since the root node of the tree, indexed by say $ij$, also cannot be color-branching, one of $i$ or $j$ must be free and thus, $b_{ij}$ can be expressed as either $b^{\{1\}}_i$ or $b^{\{2\}}_j$ in the order condition.

Conversely, we ask: given an NPRK method $(a,b,c)$ whose underlying coefficients $(a^{\{1\}},b^{\{1\}},c)$, $(a^{\{2\}},b^{\{2\}},c)$ satisfy the ARK order conditions to order $p$, does the NPRK method satisfy the NPRK order conditions to order $p$, for a generally nonlinear partition? For $p \leq 2$, the answer is affirmative since there are no color-branching trees of order $2$. For $p \geq 3$, the answer is negative. This can be seen from the color-branching trees, since we know the order conditions for the non-color-branching trees will be satisfied. It suffices to consider the third-order condition 
$$\sum_{ijkluv} b_{ij}a_{ikl}a_{juv} = 1/3$$
corresponding to the color-branching tree 
$$ \begin{tikzpicture}[baseline={0}]
    \filldraw[black] (-0.4,0.4) circle (2pt);
    \filldraw[black] (0,-0.4) circle (2pt);
    \filldraw[black] (0.4,0.4) circle (2pt);
    \draw[black, dashed, ultra thick] (-0.4,0.4) -- (0,-0.4);
    \draw[black, ultra thick] (0.4,0.4) -- (0,-0.4);
\end{tikzpicture}. $$
Let $(a^{\{1\}},b^{\{1\}},c)$ and $(a^{\{2\}},b^{\{2\}},c)$ be coefficients for a (at least) third-order ARK method. Define $a_{ijk}$ and $b_{ij}$ by \eqref{eq:ark-to-nprk-coeff-a}-\eqref{eq:ark-to-nprk-coeff-b}. Substituting \eqref{eq:ark-to-nprk-coeff-a}-\eqref{eq:ark-to-nprk-coeff-b} into the above order condition yields
\begin{align*}
\frac{1}{3} &= \sum_{ijkluv} b_{ij}a_{ikl}a_{juv} = \sum_{ijkluv} \left(\frac{b^{\{1\}}_i}{s} + \frac{b^{\{2\}}_j}{s} - \frac{c_i}{s^2}\right) a_{ikl}a_{juv} \\ 
&= \frac{1}{s} \sum_{ij} b^{\{1\}}_i c_i c_j + \frac{1}{s} \sum_{ij} b^{\{2\}}_j c_j c_i - \frac{1}{s^2} \sum_{ij}c_ic_j \\
&= \frac{1}{2s} \sum_j c_j + \frac{1}{2s} \sum_i c_i - \frac{1}{s^2} \left(\sum_i c_i\right)^2,
\end{align*}
where in the second line we used the second-order conditions $\sum_ib^{\{1\}}_ic_i = 1/2 = \sum_jb^{\{2\}}_jc_j.$ Letting $x := \sum_i c_i$, the order condition can be expressed
$$ \frac{1}{3} = \frac{x}{s} - \frac{x^2}{s^2}. $$
Viewed as a quadratic in $x$, the discriminant is $-1/(3s^2)$ which is negative since $s \in \mathbb{R} \setminus \{0\}$. Thus, the above has no real solutions, so the condition cannot be satisfied for any choice of $x$ since $c_i \in \mathbb{R}$.

Thus, we have shown that the underlying methods of an NPRK method satisfying ARK order conditions to order $p$ is not sufficient for the NPRK method to have order $p$, for $p \geq 3$. In essence, this is due to the nonlinear coupling in an NPRK method with a generally nonlinear partition $F(Y_1,Y_2)$ which creates nonlinear order conditions corresponding to color-branching trees, starting at order $p \geq 3$; for ARK methods, the partition has a linearly separable structure which removes such color-branching trees from consideration since their elementary differentials vanish.

We cannot thus in general use order $p$ ARK methods to construct order $p$ NPRK methods via equations \eqref{eq:ark-to-nprk-coeff-a}-\eqref{eq:ark-to-nprk-coeff-b}. However, since an ARK method does not arise from a unique NPRK method, it is still possible that some order $p$ ARK methods can produce an order $p$ NPRK method. 

Specifically, for order $3$, consider an order $3$ ARK method with coefficients 
$$(a^{\{1\}},b^{\{1\}},c), (a^{\{2\}},b^{\{2\}},c) \text{ such that } b^{\{1\}} = b^{\{2\}},$$ 
such as a third-order Lobatto IIIA-IIIB pair. Then, we can define an NPRK method with coefficients $a$ via equation \eqref{eq:ark-to-nprk-coeff-a} and $b$ via equation \eqref{eq:ark-to-nprk-coeff-diagonal}. It is straightforward to check that this method does satisfy the third-order condition for the order 3 color-branching tree
$$\sum_{ijkluv} b_{ij}a_{ikl}a_{juv} = 1/3,$$
since diagonality of $b$ allows one to write this condition in terms of non-color-branching order conditions.

For orders higher than $3$, using a diagonal $b$ will not work in general, due to color-branching trees such as
$$ \tau = \begin{tikzpicture}[baseline={0}]
    \filldraw[black] (-0.4,0.8) circle (2pt);
    \filldraw[black] (0,0) circle (2pt);
    \filldraw[black] (0.4,0.8) circle (2pt);
    \filldraw[black] (0,-0.8) circle (2pt);
    \draw[black, dashed, ultra thick] (-0.4,0.8) -- (0,0);
    \draw[black, ultra thick] (0.4,0.8) -- (0,0);
    \draw[black, dashed, ultra thick] (0,0) -- (0,-0.8);
\end{tikzpicture},$$
(i.e., color-branching trees which contains a non-rooted color-branching node) whose elementary weight is given by
$$ \phi_{ij}(\tau) = \sum_{klmnuv} a_{ikl}a_{kmn}a_{luv}. $$
Here, $a_{ikl}$ can not be expressed in terms of either $a^{\{1\}}$ or $a^{\{2\}}$ since $k$ and $l$ are not free. Despite the fact that an NPRK method does not necessarily have order $p$ even if its underlying integrators have order $p$ due to new nonlinear order conditions, the increased dimensionality of the NPRK tableaux and the reduced number of order conditions, compared to ARK methods, could lead to more flexibility in constructing higher-order NPRK methods. We will explore such constructions in subsequent work.

The generalization for relating NPRK$_M$ order conditions to $\text{ARK}_M$ order conditions follows similarly. Namely, the order conditions can be decomposed into the underlying RK order conditions for trees of a single color, linearly separable coupling conditions corresponding to trees which have more than one edge color but are non-color-branching, and nonlinearly separable coupling order conditions corresponding to color-branching trees. This again follows from the fact that the elementary differential for a color-branching tree vanishes, where now edges can have one of $M$ colors, assuming an $M-$additive partition 
$$ F(Y_1,\dots,Y_M) = \sum_{i=1}^M f_i(Y_i).$$
An immediate corollary of this discussion is that an NPRK$_M$ method (with generally nonlinear partition) is second-order if its $M$ underlying RK methods are, since there are no trees of order $2$ with more than one edge color.

\newpage
\section{Numerical example}\label{sec:example}
In our companion paper \cite{nprk1}, we demonstrated the utility of NPRK methods for solving the viscous Burgers' and gray thermal radiation transport equations, particularly with regard to numerical stability. We presented several NPRK methods with differing levels of stability, order, and implicitness. Here, we will instead consider an example which exemplifies the nonlinear order condition theory developed in this paper. The code for this numerical example is available at \cite{githubNPRK2025}.

\begin{figure}[H]
\begin{center}
\includegraphics[width=100mm]{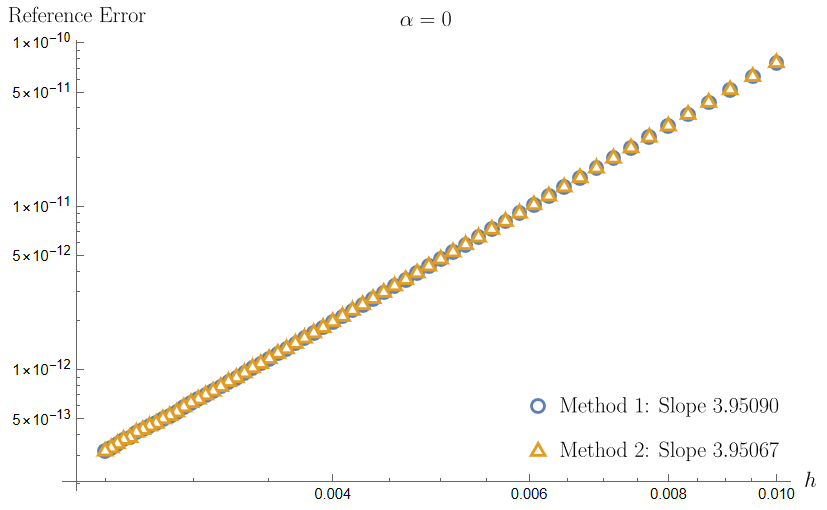}
\includegraphics[width=100mm]{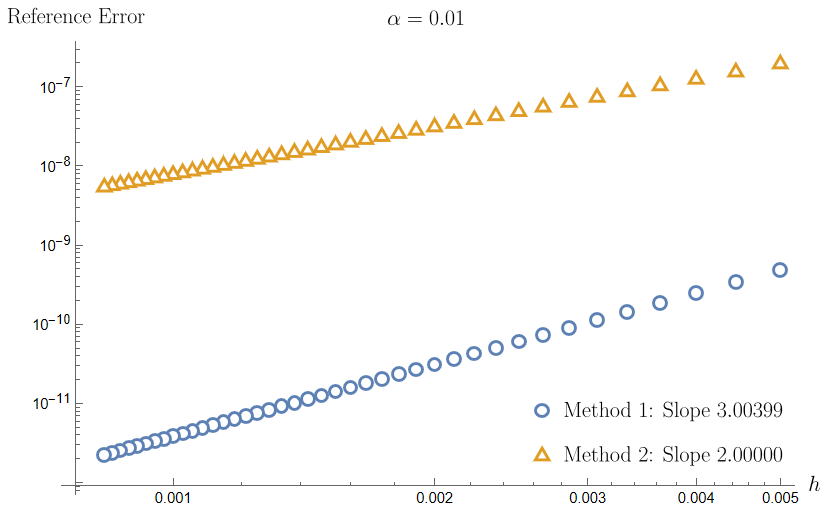}
\includegraphics[width=100mm]{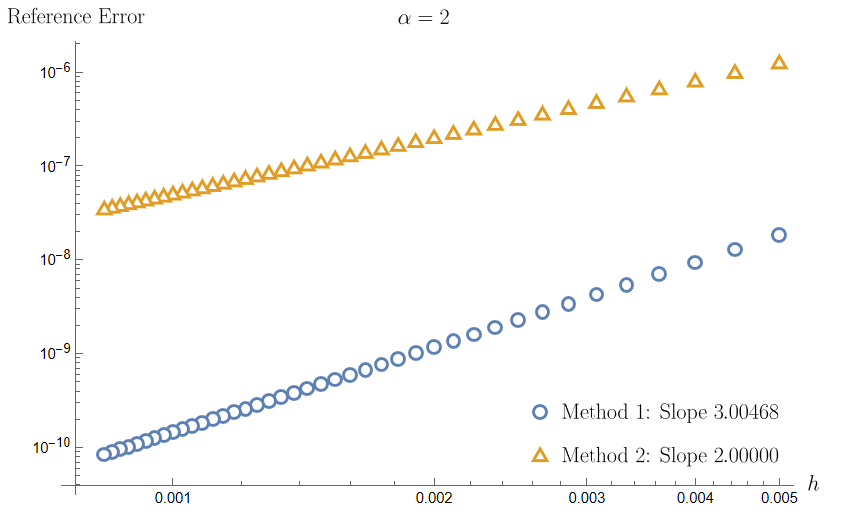}
\caption{Convergence test for Method 1 and Method 2 for $\alpha = 0, 0.01, 2$ with $u(0)=1=v(0)$ and final time $T=1$. The reference solution was computed using Method 1 with $h = 10^{-4}$. The slopes were determined using linear regression.} 
\label{figure:conv-plot}
\end{center}
\end{figure}

We will show how the nonlinear order condition theory, discussed in \Cref{section:relation-to-ark-order-cond}, can be utilized to measure the state-dependent nonlinear coupling strength in a system. We consider the Lotka--Volterra system, with $y := (u,v)^T \in \mathbb{R}^2$,
\begin{align}\label{eq:Lotka-Volterra}
    \frac{d}{dt} \begin{pmatrix} u \\ v \end{pmatrix} = \begin{pmatrix}
         u - \alpha uv \\
         v + \alpha uv
    \end{pmatrix},
\end{align}
where $\alpha \in \mathbb{R}$. Consider the Lobatto IIIA-IIIB pair with $s=3$,
\begin{align*} 
    \begin{tabular}{|c }
         $a^{\{1\}}$  \\
         \hline $b^{\{1\}}$
    \end{tabular}  =
    \begin{tabular}{|c c c}
         0 & 0 & 0  \\
         5/24 & 1/3 & -1/24 \\
         1/6 & 2/3 & 1/6 \\
         \hline 1/6 & 2/3 & 1/6
    \end{tabular}, \quad     \begin{tabular}{|c }
         $a^{\{2\}}$  \\
         \hline $b^{\{2\}}$
    \end{tabular}  =
    \begin{tabular}{|c c c}
         1/6 & -1/6 & 0  \\
         1/6 & 1/3 & 0 \\
         1/6 & 5/6 & 0 \\
         \hline 1/6 & 2/3 & 1/6
    \end{tabular}.
\end{align*}
As is well known, this method has additive order four (Section II.2.2, Theorem 2.2, of \cite{Hairer.2002}). We use this additive method to construct two NPRK methods as follows. We define Method 1 with NPRK tableau $a_{ijk}$ given by \eqref{eq:ark-to-nprk-coeff-a} and $b_{ij}$ given by \eqref{eq:ark-to-nprk-coeff-diagonal}. We define Method 2 with NPRK tableau $a_{ijk}$ again given by \eqref{eq:ark-to-nprk-coeff-a} and $\widetilde{b}_{ij}$ now given by \eqref{eq:ark-to-nprk-coeff-b}. Method 1 and Method 2 together can be interpreted as an \textit{embedded} NPRK method, since they share the same $a_{ijk}$ but have differing $b_{ij}$. As discussed in the previous section, both methods reduce to the same fourth-order additive method when the partition is additive. For a nonlinear partition, Method 1 will generally be a third-order NPRK method and Method 2 will generally be a second-order NPRK method. First, we verify the order of both methods with the above Lotka--Volterra system \eqref{eq:Lotka-Volterra}, with nonlinear partition
\begin{equation}\label{eq:LV-nonlinear-partition}
     F\left(\begin{matrix} u_1 \\ v_1 \end{matrix},\ \begin{matrix} u_2 \\ v_2  \end{matrix}\right) = \begin{pmatrix}
    u_2 - \alpha u_1 v_2 \\
    v_1 + \alpha u_2 v_1
\end{pmatrix}.
\end{equation}
See \Cref{figure:conv-plot}. As expected, when $\alpha = 0$, the partition \eqref{eq:LV-nonlinear-partition} becomes additive and thus, both methods are equivalent and exhibit fourth order convergence. When $\alpha \neq 0$, Method 1 exhibits third order convergence, whereas Method 2 exhibits second order convergence.
\newpage

\begin{figure}[H]
\begin{center}
\includegraphics[width=120mm]{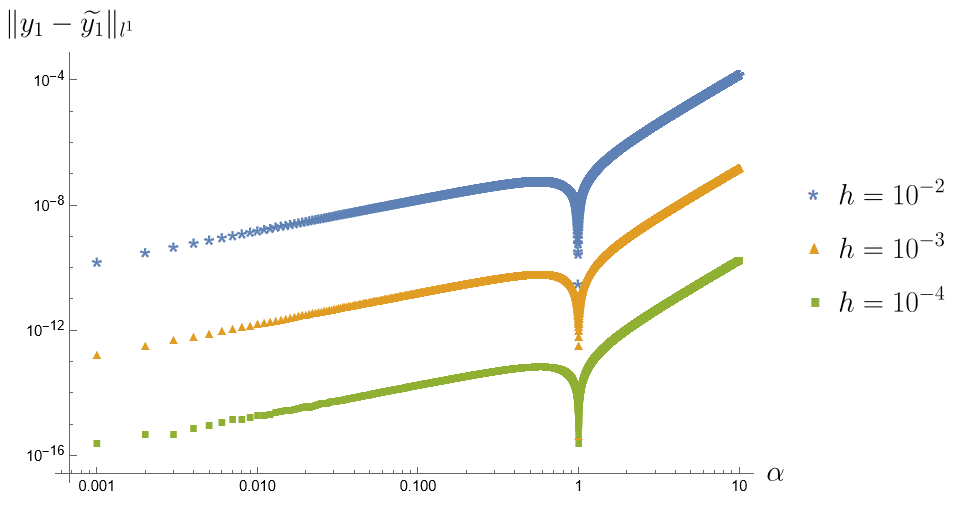}
\caption{Embedded estimate \eqref{eq:nonlinear-coupling-strength} of nonlinear coupling after one step as a function of $\alpha$ for several choices of step size $h = 10^{-2}, 10^{-3}, 10^{-4}$.} 
\label{figure:scale-test}
\end{center}
\end{figure}

Now, since both methods reduce to the same additive method, this means that the leading order difference in their accuracy is due to nonlinear coupling. Particularly, denoting $y_{n+1}$ as the numerical solution for Method 1 and $\widetilde{y}_{n+1}$ as the numerical solution for Method 2, both updated from the same $y_n$, we have an embedded estimate of the nonlinear coupling of the form
\begin{equation}\label{eq:nonlinear-coupling-strength} \|y_{n+1} - \widetilde{y}_{n+1}\| \leq C h^3 \left\|\mathcal{F}\left(\begin{tikzpicture}[baseline={0}]
    \filldraw[black] (-0.4,0.4) circle (2pt);
    \filldraw[black] (0,-0.4) circle (2pt);
    \filldraw[black] (0.4,0.4) circle (2pt);
    \draw[black, dashed, ultra thick] (-0.4,0.4) -- (0,-0.4);
    \draw[black, ultra thick] (0.4,0.4) -- (0,-0.4);
\end{tikzpicture} \right)\right\| + \mathcal{O}(h^4). 
\end{equation}
To demonstrate this embedded estimate, we compute one step of the update by Method 1 and Method 2 for the Lotka--Volterra system with $u(0)=1=v(0)$ as a function of $\alpha$ for several choices of step size $h$. The difference in the $l^1$ norm of the solutions $\|y_1-\widetilde{y}_1\|_{l^1}$ as a function of $\alpha$ is shown in \Cref{figure:scale-test}.
As expected, for each choice of $h$, as $\alpha \rightarrow 0$, we have $\|y_1 - \widetilde{y}_1\|_{l^1} \rightarrow 0$ since in this limit, $x$ and $y$ decouple and the partition becomes additive. From here, one might expect that the estimate $\|y_1 - \widetilde{y}_1\|_{l^1}$ increases monotonically as $\alpha$ increases. However, this is not the case as seen in \Cref{figure:scale-test}. Particularly, the nonlinear coupling strength also approaches zero as $\alpha \rightarrow 1$. To see why this is, note that at the initial state $u(0)=1=v(0)$, for $\alpha = 1$, $du/dt = 0$ and thus, $dv/dt = v + \alpha u v$ is no longer nonlinearly coupled to the dynamics of $u$, since $u$ is constant. In other words, the nonlinear coupling strength in a system is state-dependent, which can be measured using the embedded estimate, \eqref{eq:nonlinear-coupling-strength}.

\begin{figure}[H]
\begin{center}
\includegraphics[width=130mm]{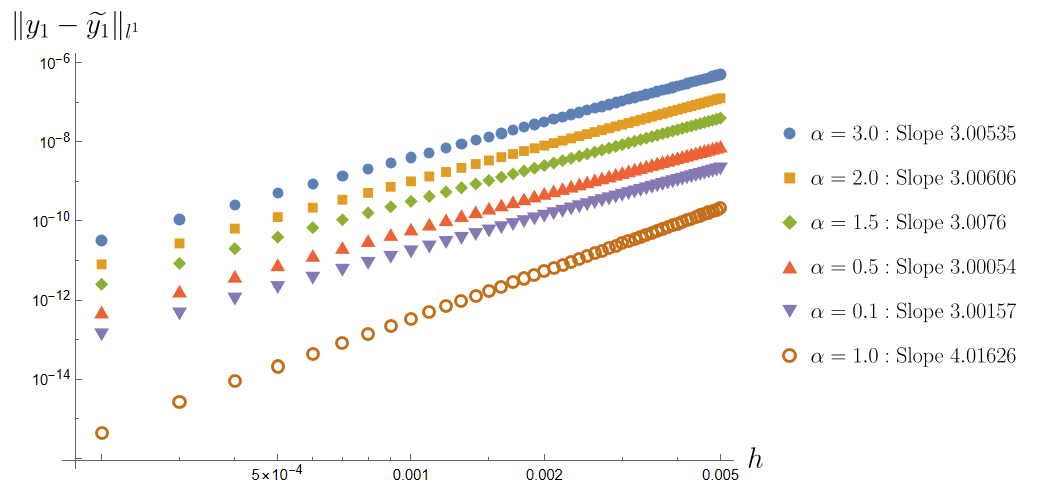}
\caption{Asymptotic scaling of the embedded estimate \eqref{eq:nonlinear-coupling-strength} with respect to $h$ for various choices of $\alpha$. The slopes were determined using linear regression.} 
\label{figure:h-scale-test}
\end{center}
\end{figure}

Stated another way, the leading term in \eqref{eq:nonlinear-coupling-strength} vanishes at the state $u(0)=1=v(0)$ as $\alpha \rightarrow 1$. To confirm this, we numerically verify the asymptotic scaling of the embedded estimate \eqref{eq:nonlinear-coupling-strength} as a function of $h$, for various choices of $\alpha$. We compute $\|y_1 - \widetilde{y}_1\|_{l^1}$ for $u(0)=1=v(0)$ as a function of $h$, where $y_1$ and $\widetilde{y}_1$ are again the solution after one step of size $h$ from Method 1 and Method 2, respectively. This is done for the choices $\alpha = 0.1, 0.5, 1.0, 1.5, 2.0, 3.0$, shown in \Cref{figure:h-scale-test}.

As expected from \eqref{eq:nonlinear-coupling-strength}, for $\alpha = 0.1, 0.5, 1.5, 2.0, 3.0$, the embedded estimate $\|y_1 - \widetilde{y}_1\|_{l^1}$ exhibits third order convergence in $h$. Interestingly, for $\alpha = 1.0$, the embedded estimate exhibits fourth order convergence. This is again explained by the fact that the leading order term of the embedded estimate \eqref{eq:nonlinear-coupling-strength} measures the nonlinear coupling which vanishes at the current state $u(0)=1=v(0)$ for $\alpha = 1$. 
Thus, in the case $\alpha = 1$ with $u(0)=1=v(0)$, \eqref{eq:nonlinear-coupling-strength} reduces to $\|y_{n+1} - \widetilde{y}_{n+1}\| \leq \mathcal{O}(h^4).$

\section*{Conclusion}
In many realistic multiphysics applications, there is nonlinear coupling of scales and physics (e.g., see \cite{southworth2023implicit,imex-trt}), making the resulting equations not directly amenable to classical additive partitions. In this paper, we provide a complete analysis of order conditions for the newly proposed class of NPRK methods via edge-colored rooted trees, which facilitates high-order partitioned integration for nonlinearly partitioned equations. General order conditions are provided for arbitrary order and number of partitions, and explicit conditions on the NPRK tableaux are provided for up to 4th-order for two partitions, and 3rd order for three partitions. NPRK order conditions are also related to ARK order conditions, and it is shown that for order $p\geq 3$, additional conditions related to nonlinear coupling cannot be represented in trees used to derive ARK order conditions (although it is possible for an ARK method of order three to satisfy the nonlinear coupling conditions when posed as an NPRK method). 

Here we provide examples of a second-order and third-order NPRK method constructed from the additive Lobatto IIIA-IIIB pair, differing only in the nonlinear order conditions, and show how they can be used to obtain an estimate of state-dependent nonlinear coupling strength. Our companion paper provides a number of other example methods of second and third order \cite{nprk1}. In future work, we will derive practical NPRK methods of various orders, structures, and number of partitions, focusing on optimizing coefficients for properties related to accuracy and stability. We will also explore the NPRK tensor structure as a simplifying assumption for deriving ARK methods of higher order and number of partitions, due to the reduced number of order conditions compared with classical ARK methods. We have not addressed structure preservation here; future work will also study method properties such as conservation, symplecticity, adaptivity.

\section*{Acknowledgements}
The authors would like to thank the reviewers for their helpful comments and suggestions. BKT was supported by the Marc Kac Postdoctoral Fellowship at the Center for Nonlinear Studies at Los Alamos National Laboratory. BSS was supported by the Laboratory Directed Research and Development program of Los Alamos National Laboratory under project number 20220174ER. Los Alamos National Laboratory report LA-UR-24-20566. TB was funded by the National Science Foundation under grant NSF-OIA-2327484.

\begin{appendix}

\section{List of Order Conditions}\label{app:order-list}
Below are lists of NPRK$_M$ order conditions for $M=2$, third-order \eqref{eq:n2order3list} and fourth-order \eqref{eq:n2order4list}, and $M=3$ third-order \eqref{eq:n3order3list}, generated by a symbolic version of \Cref{alg:two-nprk-order}. The summations are understood to run over all present indices from $1$ to $s$. Equations annotated with ${}^*$ denote linearly separable coupling conditions and equations annotated with ${}^\dagger$ denote nonlinear coupling conditions. The code for all results generated in this paper is available at \cite{githubNPRK2025}.

\begin{equation}\label{eq:n2order3list}
\text{Third-order conditions for } M=2
\end{equation}
\begin{equation*}
\begin{split}
\sum b_{ij}a_{ikl} a_{kuv}  &=\frac{1}{6},\quad  \\
{}^*\sum b_{ij}a_{kuv} a_{jkl}  &=\frac{1}{6},\quad  \\
{}^*\sum b_{ij}a_{ikl} a_{luv}  &=\frac{1}{6},\quad  \\
\sum b_{ij}a_{jkl} a_{luv}  &=\frac{1}{6},\quad  
\end{split} 
\begin{split}
\sum b_{ij}a_{ikl} a_{iuv}  =\frac{1}{3},\\
{}^\dagger \sum b_{ij}a_{iuv} a_{jkl}  =\frac{1}{3}, \\
\sum b_{ij}a_{jkl} a_{juv}  =\frac{1}{3}.
\end{split}
\end{equation*}

\newpage
\begin{equation}\label{eq:n2order4list}
\text{Fourth-order conditions for } M=2
\end{equation}
\begin{equation*}
\begin{split}
\sum b_{i j} a_{i k l} a_{k u v} a_{u a b}  &= \frac{1}{24},\quad  \\
{}^*\sum b_{i j} a_{k u v} a_{u a b} a_{j k l}  &= \frac{1}{24},\quad  \\
{}^*\sum b_{i j} a_{i k l} a_{u a b} a_{l u v}  &= \frac{1}{24},\quad  \\
{}^*\sum b_{i j}a_{u a b} a_{j k l} a_{l u v} &= \frac{1}{24},\quad  \\
{}^*\sum b_{i j}a_{i k l} a_{k u v} a_{v a b} &=\frac{1}{24},\quad  \\
{}^*\sum b_{i j}a_{k u v} a_{j k l} a_{v a b} &=\frac{1}{24},\quad  \\
{}^*\sum b_{i j}a_{i k l} a_{l u v} a_{v a b} &=\frac{1}{24},\quad  \\
\sum b_{i j}a_{j k l} a_{l u v} a_{v a b} &=\frac{1}{24},\quad  \\
\sum b_{i j}a_{i k l} a_{k u v} a_{k a b} &=\frac{1}{12},\quad  \\
{}^*\sum b_{i j}a_{k u v} a_{k a b} a_{j k l} &=\frac{1}{12},\quad  \\
{}^\dagger\sum b_{i j}a_{i k l} a_{k a b} a_{l u v} &=\frac{1}{12},\quad  \\
{}^\dagger\sum b_{i j}a_{k a b} a_{j k l} a_{l u v} &=\frac{1}{12},\quad  \\
{}^*\sum b_{i j}a_{i k l} a_{l u v} a_{l a b} &=\frac{1}{12},\quad  
\end{split}
\begin{split}
\sum b_{i j}a_{j k l} a_{l u v} a_{l a b} &=\frac{1}{12},  \\
\sum b_{i j}a_{i k l} a_{i a b} a_{k u v} &=\frac{1}{8},  \\
{}^\dagger\sum b_{i j}a_{i a b} a_{k u v} a_{j k l} &=\frac{1}{8},  \\
{}^*\sum b_{i j}a_{i k l} a_{i a b} a_{l u v} &=\frac{1}{8}, \\
{}^\dagger\sum b_{i j}a_{i a b} a_{j k l} a_{l u v} &=\frac{1}{8},  \\
{}^\dagger\sum b_{i j}a_{i k l} a_{k u v} a_{j a b} &=\frac{1}{8},  \\
{}^*\sum b_{i j}a_{k u v} a_{j k l} a_{j a b} &=\frac{1}{8},  \\
{}^\dagger\sum b_{i j}a_{i k l} a_{j a b} a_{l u v} &=\frac{1}{8},  \\
\sum b_{i j}a_{j k l} a_{j a b} a_{l u v} &=\frac{1}{8},  \\
\sum b_{i j}a_{i k l} a_{i u v} a_{i a b} &=\frac{1}{4},  \\
{}^\dagger\sum b_{i j}a_{i u v} a_{i a b} a_{j k l} &=\frac{1}{4},  \\
{}^\dagger\sum b_{i j}a_{i a b} a_{j k l} a_{j u v} &=\frac{1}{4},  \\
\sum b_{i j}a_{j k l} a_{j u v} a_{j a b} &=\frac{1}{4}.
\end{split}
\end{equation*}

\begin{equation}\label{eq:n3order3list}
\text{Third-order conditions for } M=3
\end{equation}
\begin{equation*}
\begin{split}
\sum b_{i j k}a_{i u v w} a_{u a b c}  &= \frac{1}{6},\quad  \\
{}^*\sum b_{i j k}a_{i u v w} a_{v a b c} &=\frac{1}{6},\quad  \\
{}^*\sum b_{i j k}a_{i u v w} a_{w a b c} &=\frac{1}{6},\quad  \\
{}^*\sum b_{i j k}a_{u a b c} a_{j u v w} &=\frac{1}{6},\quad  \\
\sum b_{i j k}a_{j u v w} a_{v a b c} &=\frac{1}{6},\quad  \\
{}^*\sum b_{i j k}a_{j u v w} a_{w a b c} &=\frac{1}{6},\quad  \\
{}^*\sum b_{i j k}a_{u a b c} a_{k u v w} &=\frac{1}{6},\quad  \\
{}^*\sum b_{i j k}a_{v a b c} a_{k u v w} &=\frac{1}{6},\quad  
\end{split}
\begin{split}
\sum b_{i j k}a_{k u v w} a_{w a b c} &=\frac{1}{6},  \\
\sum b_{i j k}a_{i u v w} a_{i a b c} &=\frac{1}{3},  \\
{}^\dagger\sum b_{i j k}a_{i a b c} a_{j u v w} &=\frac{1}{3},  \\
{}^\dagger\sum b_{i j k}a_{i a b c} a_{k u v w} &=\frac{1}{3},  \\
\sum b_{i j k}a_{j u v w} a_{j a b c} &=\frac{1}{3},  \\
{}^\dagger\sum b_{i j k}a_{j a b c} a_{k u v w} &=\frac{1}{3},  \\
\sum b_{i j k}a_{k u v w} a_{k a b c} &=\frac{1}{3}.
\end{split}
\end{equation*}

\end{appendix}

\newpage
\bibliographystyle{plainnat}
\bibliography{refs}

\end{document}